\documentclass[11pt]{amsart}
\usepackage[utf8]{inputenc}
\usepackage[english]{babel}

\usepackage{amscd,amssymb}
\usepackage{amsmath}
\usepackage{amsthm}
\usepackage{graphicx}
\usepackage{fancybox}
\usepackage{epic,eepic}
\usepackage{amstext}
\usepackage{mathtools}
\usepackage{esint}
\usepackage[square,numbers]{natbib}
\usepackage{booktabs}
\usepackage{hyperref}
\usepackage[capitalise,nameinlink]{cleveref}
\crefformat{equation}{\textup{#2(#1)#3}}
\usepackage{comment}
\usepackage{mathtools}
\usepackage{tikz,pgfplots}
\usepackage{subcaption}

\def\K{\mathcal{K}}

\def\H{\mathcal{H}}

\def\R{\mathcal{R}}

\def\K{\mathcal{K}}

\def\Pnull{\mathbb{P}^0}

\def\elem{T}

\def\T{\mathcal{T}}
\def\TH{\mathcal{T}_H}

\def\PiH {\Pi_H}

\def\linh{\mathrm{span}}

\def\one{\mathbf{1}}

\newcommand{\mnormf}[1]{{| #1 |}_{M}}
\newcommand{\mnormfo}[2]{{| #1 |}_{M,#2}}

\newcommand{\lnormf}[1]{{| #1 |}_{L}}
\newcommand{\lnormfo}[2]{{| #1 |}_{L,#2}}

\newcommand{\knormf}[1]{{| #1 |}_{K}}

\newcommand{\vsp}[2]{{\left( #1\,,\,#2 \right)}}
\newcommand{\vspf}[2]{{( #1\,,\,#2 )}}
\newcommand{\vnormo}[2]{{\left\vert #1 \right\vert}_{V,#2}}
\newcommand{\vnormfo}[2]{{| #1 |}_{V,#2}}

\newcommand{\vnormf}[1]{{| #1 |}_{V}}

\newcommand{\with}{\,:\,}
\DeclareMathOperator*{\argmin}{argmin}

\definecolor{myBlue}{RGB}{113,104,238} 
\definecolor{myGreen}{RGB}{114,175,30} 
\definecolor{myRed}{RGB}{180,50,50}  
\definecolor{myOrange}{RGB}{225,92,22}

\newtheorem{theorem}{Theorem}[section]

\newtheorem{lemma}[theorem]{Lemma}
\newtheorem{assumption}[theorem]{Assumption}

\theoremstyle{definition}

\theoremstyle{remark}
\newtheorem{remark}[theorem]{Remark}
\numberwithin{theorem}{section}
\numberwithin{equation}{section}
\numberwithin{table}{section}
\numberwithin{figure}{section}

\allowdisplaybreaks

\begin{document}
	
\title[]{Super-localization of spatial network models}
\author[]{Moritz Hauck$^\dagger$, Axel M\aa lqvist$^{\ddagger}$}
\address{${}^{\dagger}$ Institute of Mathematics, University of Augsburg, Universit\"atsstr.~12a, 86159 Augsburg, Germany}
\email{moritz.hauck@uni-a.de}
\address{${}^{\ddagger}$ Department of Mathematical Sciences, Chalmers University of Technology and University of Gothenburg, 41296 Göteborg, Sweden}
\email{axel@chalmers.se}
\thanks{The work of Moritz Hauck is part of a project that has received funding from the European Research Council (ERC) under the European Union's Horizon 2020 research and innovation programme (Grant agreement No.~865751 --  RandomMultiScales). The work of Axel M\aa lqvist is supported by the Swedish research council (project number 2019-03517).}
\maketitle

\begin{abstract}
	
Spatial network models are used as a simplified discrete representation in a wide range of applications, e.g., flow in blood vessels, elasticity of fiber based materials, and pore network models of porous materials. Nevertheless, the resulting linear systems are typically large and poorly conditioned and their numerical solution is challenging.  

This paper proposes a numerical homogenization technique for spatial network models which is based on the Super Localized Orthogonal Decomposition (SLOD), recently introduced for elliptic multiscale partial differential equations. It provides accurate coarse solution spaces with approximation properties 
independent of the smoothness of the material data. A unique selling 
point of the SLOD is that it constructs an almost local basis of these coarse spaces, requiring less computations on the fine scale and achieving improved sparsity on the coarse scale compared to other state-of-the-art methods. 
We provide an a-posteriori analysis of the proposed method and numerically confirm the method's unique localization properties. In addition, we show its applicability also for high-contrast channeled material data.
\end{abstract}

\vspace{1cm}
\\[2ex]

\section{Introduction}

Spatial networks are a useful tool for constructing simplified discrete representations of complex geometric structures. Blood vessels may, for example, be modeled as connected one dimensional line segments forming a network of nodes and edges, see \cite{FKOWW22}. Paper consists of a web of wooden fibers that may be modeled as one dimensional beams forming a network, see \cite{KMM20}. Also  porous materials, such as sandstone, may be represented by a pore network model, see \cite{EPT12}. Such simplifications reduce the complexity from a full three dimensional geometry to a discrete model for which computer simulations can be performed more easily. Nevertheless, highly heterogeneous materials and complicated geometries typically cause the underlying linear systems to be poorly conditioned. 

This paper considers spatial network models that can be described by weighted graph Laplacians, arising from applications modeled by elliptic partial differential equations (PDEs). There already exists a variety of well-established iterative multilevel solvers for such problems.  A prominent example is algebraic multigrid, see e.g.~\cite{LB12,XZ17,HWZ21} and the references therein. In a purely algebraic setting, the construction of multiple discretization levels can be challenging. For spatial networks with sufficient structure, it is possible to embed the network into a domain $\Omega \subset \mathbb R^d$ and introduce scales by interpolating the network onto a family of finite element meshes, see e.g.~\cite{GoHeMa22}. On coarser scales, the spatial network behaves like a continuous material  and therefore inherits many advantageous properties from the continuous setting. This can be utilized for transferring successful algorithms for solving PDEs to spatial network models. We will focus on numerical homogenization, with the goal to compute an accurate effective coarse scale model of the full network.

Numerical homogenization is about the construction of problem-adapted, optimally approximating ansatz spaces  possessing almost local computable bases. Near optimal numerical homogenization is achieved by the Localized Orthogonal Decomposition (LOD) \cite{MaP14,HeP13,KPY18,MalP20,BrennerLOD,Peterseim2021} and Gamblets \cite{Owh17,OwhS19} which construct a fixed number of basis functions per mesh entity that decay exponentially relative to the mesh. For the computation of the basis, a localization of the basis functions to element patches is performed. The number of element layers in the element patches needs to be increased logarithmically with the desired accuracy. 
An alternative approach is taken by the AL-basis \cite{GGS12} and G(ms)FEM methods \cite{BaL11,Ma21} which solve local spectral problems  and construct the global ansatz space using a partition of unity. Here, the support of the basis functions is fixed by the choice of partition of unity and the number of local eigenfunctions needs to be increased logarithmically with the desired accuracy. 

Recently, the Super Localized Orthogonal Decomposition (SLOD) has been proposed in \cite{HaPe21b} performing a localization of the same space as the LOD, but utilizing a novel localization strategy allowing for super-exponentially decaying basis functions. This improved localization results in smaller local patch problems for the basis computation and a sparser coarse system matrix. The method proved its effectiveness also beyond elliptic multiscale problems, see  \cite{FHP21,BFP22}. 
  
There are some contributions specifically targeting numerical homogenization of spatial network models, see \cite{EILRW09,ILW10,KMM20,EGHKM22}. In particular, a spatial network version of the LOD has been developed in \cite{EGHKM22}. Therein a rigorous proof of the exponential decay of the LOD basis functions in the spatial network setting is provided utilizing the techniques developed in \cite{KorY16} and \cite{GoHeMa22}.

In this paper, we develop a SLOD for spatial network models. As model problem, we consider a weighted graph Laplacian $K$ on the spatial network $\mathcal{G}=(\mathcal{N},\mathcal{E})$, i.e., we seek the solution to the possibly poorly conditioned linear system 
\begin{equation}
	\label{eq:modelproblem0}
	Ku=Mf,
\end{equation}
where $u$ fulfills homogeneous Dirichlet boundary conditions, $M$ is a diagonal mass matrix, and $f$ is a given source term. For a depiction of a possible spatial network, see \Cref{fig:spatialnetwork}. 
Let the spatial network $\mathcal G$ be embedded into a domain $\Omega$; for the ease of presentation, let $\Omega = [0,1]^d$. We consider coarse finite element meshes $\TH$ of $\Omega$ and define problem-adapted ansatz spaces by applying the solution operator $\mathcal K^{-1}$ of \eqref{eq:modelproblem0} to $\TH$-piecewise constants. Such ansatz spaces have uniform approximation rates, independently of the smoothness of the material data. The SLOD then identifies local $\TH$-piecewise constant right-hand sides that yields a rapidly decaying response under the operator $\mathcal K^{-1}$. These responses are then localized to element patches and used as basis functions of the SLOD. 

We derive an a-posteriori error bound for the SLOD error in terms of the size of the element patches and the coarse mesh size which holds under mild assumptions on the underlying network. In a series of numerical experiments, we illustrate the network assumptions, the super-locality of the SLOD basis and the method's performance in presence of high-contrast data, in particular channeled data.

The outline of the paper is as follows. In \Cref{s:mp}, we introduce the spatial network model. A coarse scale discretization together with assumptions on the underlying network is presented in \Cref{sec:network}. \Cref{s:idealmethod} then constructs a prototypical problem-adapted ansatz space. We identify rapidly decaying basis functions and perform a localization of the basis in \Cref{s:locapprox}. \Cref{s:locmethod} states the SLOD method along with an a-posteriori error estimate and, finally, \Cref{s:ne} presents  numerical experiments. 

\section{A spatial network model}\label{s:mp}

\begin{figure}[h]
	\includegraphics[width=.6\linewidth]{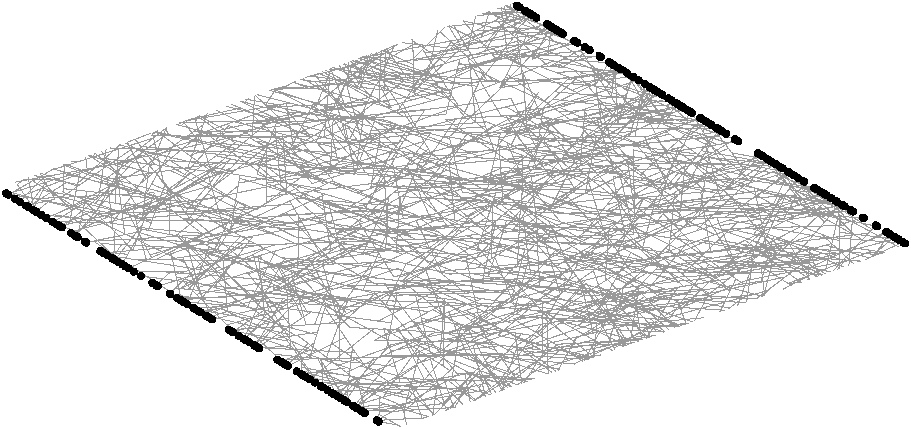}
	\caption{Illustration of spatial network with Dirichlet nodes in $\Gamma$ marked black.}
	\label{fig:spatialnetwork}
\end{figure}

We consider a connected graph $\mathcal{G} = (\mathcal{N}, \mathcal{E})$ of nodes
$\mathcal{N}$ and  edges
\begin{equation*}
\mathcal{E} = \{ \{x, y\} \,:\, \text{an edge connects nodes } x,y\in\mathcal N\}
\end{equation*}
that is embedded into the unit hyper-cube $\Omega = [0,1]^d$, $d \in \mathbb N$.
The presented methodology can also be generalized to polygonal and polyhedral domains, however, for the ease of presentation, we only consider hyper-cubes. We write $x \sim y$ if two nodes $x,y\in\mathcal{N}$ are connected by an edge in $\mathcal E$ and denote by $|x-y|$ the Euclidean distance between the nodes. By $\mathcal{N}(\omega)$, we denote all nodes that are contained in the subset $\omega\subset \Omega$. We impose homogeneous Dirichlet boundary conditions for nodes on the boundary segment  $\Gamma \subset \partial \Omega$. For a depiction of an example of a spatial network, see \Cref{fig:spatialnetwork}.

For the subsequent presentation, we introduce function spaces on the network. Let $\hat V$ denote the space of all real valued functions  defined on the set of nodes $\mathcal{N}$ and denote by
$$V = \{ v \in \hat V\with v(x) = 0,\; x \in \Gamma\}$$
the subset of $\hat {V}$ satisfying homogeneous Dirichlet boundary conditions on the boundary segment $\Gamma$. Furthermore, let $V_\omega$ and $\hat V_\omega$ denote the spaces of functions in $V$ and $\hat V$, respectively, that are supported in the subdomain $\omega$.

\subsection{Linear operators}
We define a diagonal edge length weighted mass operator $M\colon \hat V \to \hat V$ by the following sum of node contributions
\begin{align}
\label{eq:M}
M\coloneqq \sum_{x\in\mathcal{N}}M_x\quad \text{ with }\quad M_x\colon \hat V \to \hat V,\quad (M_x v, v) \coloneqq \frac{1}{2}\sum_{y \sim x} |x - y| v(x)^2.
\end{align}
It shall be noted that the node contributions $M_x$ are uniquely defined by  their associated quadratic forms. For subdomains $\omega \subset \Omega$, we define a local version of $M$ as $$M_\omega \coloneqq \sum_{x \in \mathcal N(\omega)} M_x.$$
The bilinear form $\vspf{M \cdot}{\cdot}$ is an inner product on the space $\hat V$ with induced norm $\mnormf{\cdot}^2 \coloneqq \vsp{M \cdot}{\cdot}$. By restriction to $\omega$, we obtain the semi-norm $\mnormfo{\cdot}{\omega}^2 \coloneqq \vsp{M_\omega \cdot}{\cdot}$ which can be interpreted as the mass of the network in subdomain $\omega$.

Furthermore, we define a reciprocal edge length weighted graph Laplacian $L\colon \hat V \to \hat V$~by
\begin{equation}
\label{eq:L}
L\coloneqq \sum_{x\in\mathcal{N}}L_x\quad \text{ with  }\quad L_x\colon \hat V \to \hat V,\quad  (L_x v, v) \coloneqq  \frac{1}{2}\sum_{y\sim x} \frac{(v(x) - v(y))^2}{|x - y|},
\end{equation}
where the node contributions $L_x$ are symmetric and are again uniquely defined  by their associated quadratic forms. 
A local version of $L$ is given by 
\begin{equation}
\label{eq:locL}
L_\omega \coloneqq \sum_{x \in \mathcal{N}(\omega)} L_x.
\end{equation}
It shall be noted that $(L_\omega v)(x)$ is nonzero for nodes $x$ outside of
$\omega$ that are adjacent to nodes in $\omega$.
Since the graph $\mathcal G$ is assumed to be connected, the kernel of $L$ consists only of the globally constant functions. Hence,  $\lnormf{\cdot}^2 = (L\cdot, \cdot)$ defines a semi-norm on $\hat V$. By restriction to  $\omega$, we obtain $\lnormfo{\cdot}{\omega}^2 \coloneqq  (L_\omega
\cdot, \cdot)$.

\begin{remark}[Weighting of $M$ and $L$]
	The non-standard weightings in \eqref{eq:M} and \eqref{eq:L} are chosen such that, in one dimension, the operators $L$ and $M$ resemble the first order finite element stiffness matrix of the Poisson equation and the corresponding lumped mass matrix, respectively. This weighting enables an analysis which is similar in style to well-established PDE analysis, see e.g. \cite{KMM20}. 
\end{remark}

By combining the mass matrix and the graph Laplacian,  we obtain another inner product on $\hat V$, namely $\vsp{(L+M)\,\cdot}{\cdot}$. For its induced norm, we write $\vnormf{\cdot}^2 \coloneqq \mnormf{\cdot}^2+\lnormf{\cdot}^2$ and the restriction of the norm to a subdomain $\omega$ is defined by $\vnormfo{\cdot}{\omega}^2 \coloneqq \mnormfo{\cdot}{\omega}^2 + \lnormfo{\cdot}{\omega}^2$.
\subsection{Model problem}

For demonstrating the extension of the SLOD to the spatial network setting, this paper considers a weighted graph Laplacian as model problem. It is defined as follows
\begin{equation}
\label{eq:wgraphlaplacian}
  K\coloneqq \sum_{x\in\mathcal{N}}K_x\quad \text{ with  }\quad K_x\colon \hat V \to \hat V,\quad  (K_x v, v) \coloneqq  \frac{1}{2}\sum_{y\sim x}\gamma_{xy} \frac{(v(x) - v(y))^2}{|x - y|}
\end{equation}
with edge weights 
\begin{equation}
\label{eq:edgeweights}
0<\alpha\leq \gamma_{xy}\leq\beta< \infty
\end{equation}
 determining the edge's conductivity. For subdomains 
 $\omega$,  local versions $K_\omega$ of can be defined analogously to \eqref{eq:locL}.
 
 Given a right-hand side $f \in \hat V$, we seek $u \in V$ such that, for all $v \in V$,
\begin{equation}\label{eq:weakform}
  (K u, v) = (Mf, v).
\end{equation}

The unique solvability of this problem can be ensured under the minimal requirements that the underlying network $\mathcal G$ is connected and that there exists at least one Dirichlet node in~$\Gamma$, i.e., $\mathcal N(\Gamma) \neq \emptyset$. Subsequently, these assumptions are always assumed to hold. 

Using the bounds on the edge weights \eqref{eq:edgeweights} we obtain that, for all $v \in \hat V$,
\begin{equation}
	\label{eq:normeq}
	\alpha\vspf{Lv}{v}\leq\vspf{Kv}{v}\leq \beta \vspf{Lv}{v} 
\end{equation}
which, in particular implies that on $V$, the energy norm $|\cdot|^2 \coloneqq \vspf{K\cdot}{\cdot}$ is equivalent to the norm $\lnormf{\cdot}$.

\subsection{Global Friedrichs' inequality}
It is possible to derive a Friedrichs' inequality in the spatial network setting.

\begin{lemma}[Friedrichs' inequality]
	\label{l:friedrichs}
	There exists $C_\mathrm{fr}>0$, such that, for all $v \in V$,
	\begin{equation}
	\label{eq:friedrichs}
	\mnormf{v} \leq C_{\mathrm{fr}} \lnormf{v}.
	\end{equation}
\end{lemma}
\begin{proof}
	We consider the generalized eigenvalue problem $Lv = \lambda Mv$ posed in the space $V$ with $\lambda_1 \leq \lambda_2\leq \dots$ denoting its eigenvalues. The connectivity of the spatial network $\mathcal G$ and  $\mathcal N(\Gamma) \neq \emptyset$ ensure that $\lambda_1>0$. The first eigenvalue $\lambda_1$ can be characterized by the min-max theorem as follows
	\begin{equation*}
	\lambda_1 = \inf_{0\neq v \in V}\frac{\vsp{Lv}{v}}{\vsp{Mv}{v}}>0
	\end{equation*}
	 which can be reformulated as
	\begin{equation*}
	\mnormf{v}^2 \leq \lambda_1^{-1}\lnormf{v}^2.
	\end{equation*}
	Thus, \eqref{eq:friedrichs} follows with constant $C_\mathrm{fr} = \lambda_1^{-1/2}>0$. 
\end{proof}

	Note that \Cref{l:friedrichs} does not provide an explicit characterization of $C_\mathrm{fr}$ in terms of the properties of the underlying network. However, one can establish a connection to the first eigenvalue of the well-studied normalized graph Laplacian for which such explicit characterizations exist, see e.g.~\cite{CY95,C05}. 

\subsection{Stability estimates}
Let us also introduce a local variant of \eqref{eq:weakform} which, for a subdomain $\omega \subset \Omega$ and a given  $f\in V_\omega$, seeks $u_\omega \in V_\omega$ such that, for all $v \in V_\omega$,
\begin{equation}
\label{eq:weakformlocal}
\vsp{K_\omega u_\omega }{v} = \vsp{M_\omega f}{v}.
\end{equation}
Existence and uniqueness of the solution $u_\omega$ follow since $K_\omega$, by construction, is also a weighted graph Laplacian and $V_\omega\subset V$.
Friedrichs' inequality then allows us to show the stability of the model problem~\eqref{eq:weakform} and its local counterpart \eqref{eq:weakformlocal} with respect to $\lnormf{\cdot}$ and $\lnormfo{\cdot}{\omega}$ which are norms on $V$ and~$V_\omega$, respectively. 

\begin{lemma}[Stable solvability]
	The solution to \eqref{eq:weakform} is stable in the sense that 
	\begin{equation*}
	\lnormf{u} \leq C_\mathrm{fr}\alpha^{-1} \mnormf{f}
	\end{equation*}
	with $C_\mathrm{fr}$ from \Cref{l:friedrichs}. This also holds for the local version \eqref{eq:weakformlocal}, i.e., 
	\begin{equation*}
	\lnormfo{u_\omega}{\omega} \leq C_\mathrm{fr}\alpha^{-1} \mnormfo{f}{\omega}. 
	\end{equation*}
	\end{lemma}
	\begin{proof}
		Using the uniform lower bound \eqref{eq:edgeweights},  definition \eqref{eq:weakform} and \Cref{l:friedrichs}, we obtain
	\begin{equation*}
	\alpha\lnormf{u}^2 \leq \vsp{K u}{u} = \vsp{Mf}{u} \leq \mnormf{f} \mnormf{u} \leq C_\mathrm{fr} \mnormf{f} \lnormf{u}.
	\end{equation*}
	Dividing by $\lnormf{u}$, the global stability estimate follows. 
	
	The local stability estimate can be obtained similarly noting that $u_\omega$, in particular,  is also an element of $V$ for which we can apply \Cref{l:friedrichs}. Using that $\mnormf{u_\omega} = \mnormfo{u_\omega}{\omega}$ and $\lnormf{u_\omega} = \lnormfo{u_\omega}{\omega}$, the local stability estimate can be concluded.
\end{proof}

For later use, we define the solution operator $\K^{-1}\colon \hat V \to V,$ $f \mapsto u$ mapping a right-hand side $f$ to its corresponding solution $u$ solving \eqref{eq:weakform}. On subsets $\omega \subset \Omega$, we also define the local solution operator $\K_\omega^{-1}\colon \hat V_\omega \to V_\omega$, $f\mapsto u_\omega $ mapping a local right-hand side to the local solution $u_\omega$ satisfying \eqref{eq:weakformlocal}.

\section{Coarse scale discretization and network assumptions}\label{sec:network}

\subsection{Coarse mesh} 
The proposed method constructs its problem-adapted basis functions with respect to some coarse mesh $\TH$. For simplicity, we restrict ourselves to Cartesian meshes which we define, unlike classical textbooks on finite element theory, as
\begin{equation*}
\TH \coloneqq \{S_H(x) \with  x = (x_1, \ldots, x_d) \in \Omega \text{ and } H^{-1}x_i + 1/2 \in \mathbb{Z} \text{ for } i = 1, \dots, d\}
\end{equation*}
with elements that are neither closed nor open but are rather defined as
\begin{equation*}
S_H(x) = [x_1 - H/2, x_1 + H/2) \times \cdots \times [x_d - H/2, x_d
+ H/2).
\end{equation*}
If $x_i + H/2 = 1$, we replace $[x_i - H/2, x_i + H/2)$ by $[x_i - H/2, x_i + H/2]$. This definition ensures that the elements form a true partition of $\Omega$ meaning that any point in ${\Omega}$ is contained in exactly one element. 

\begin{figure}[h]
	\includegraphics[width=.35\linewidth]{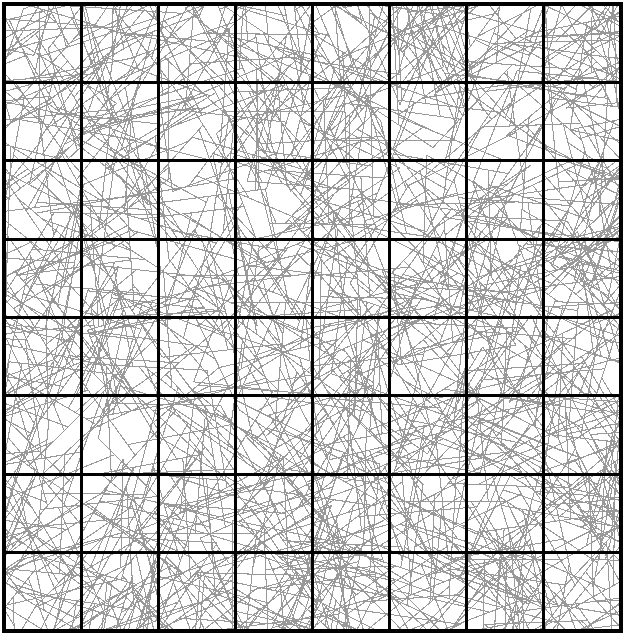}
	\caption{Coarse mesh $\T_{2^{-3}}$  with underlying spatial network.}
	\label{fig:networkwcoarsemesh}
\end{figure}

 We also introduce the concept of patches which is based on neighborhood relations of elements. The first order patch of a subset $\omega \subset \Omega$ is defined as
\begin{equation*}
\mathsf{N}(\omega) \coloneqq \{ x \in \Omega\with  \exists T \in \mathcal{T}_H \with x \in T,\;  \overline{T} \cap \overline{\omega}
\ne \emptyset\}.
\end{equation*}
By recursion, we can then define the $\ell$-th order patch as $\mathsf{N}^\ell(\omega) \coloneqq \mathsf{N}^{\ell-1}(\mathsf{N}(\omega))$ with
$\mathsf{N}^1 =\mathsf{N}$. 

\subsection{Network assumption}
The proposed method relies on the assumption that the mesh $\TH$ is coarse compared to the underlying network. On the length scale of the coarse mesh, the network then inherits many properties known from the continuous setting, e.g., an element-wise Poincar\'{e} inequality. 
 
Hence, we restrict ourselves to meshes with mesh sizes larger than some parameter $H_0>0$. More precisely, 
we consider a hierarchy of meshes 
\begin{equation}
\label{eq:hierarchy}
\{\TH \with H \in \mathcal H\}
\end{equation}
 with $\mathcal H$ denoting a finite set of positive mesh size parameters with the smallest element being $H_0$. 
The parameter $H_0$ can be interpreted as the smallest mesh size for which desired properties from the continuous case still carry over. 
The following assumption provides an insight into the choice of $H_0$.

\begin{samepage}
\begin{assumption}[Network connectivity]
  \label{ass:network}
  Assume that the smallest mesh size $H_0$ of the hierarchy of meshes \eqref{eq:hierarchy} is chosen such that, for any element $T\in \TH$, $H \in \H$, there is a connected subgraph $\bar {\mathcal G} = (\bar{\mathcal N}, \bar{\mathcal E})$ of $\mathcal G$ that contains
  \begin{itemize}
  	\item all edges with one or both endpoint in $T$,
  	\item  only edges with endpoints contained within $\mathsf N(T)$.
  \end{itemize}   
\end{assumption}
\end{samepage}

\subsection{Element-wise Poincar\'{e} inequality}
Under the previous assumption, it is possible to prove the following element-wise Poincar\'{e} inequality. 

\begin{lemma}[Element-wise Poincar\'{e}'s inequality]
	\label{l:poincare}
	Let \Cref{ass:network} be satisfied. Then, for all $T \in \TH$, $H \in \H$, there exists a constant $C_\mathrm{po}>0$ such that, for all $v \in \hat V$, there exists a constant $c$ such that
	\begin{equation}
	\label{eq:poincare}
	|v - c|_{M,T} \le C_\mathrm{po} |v|_{L,\mathsf N(T)}.
	\end{equation}
	The constant $c$ resembles the element-average in the continuous case.
\end{lemma}
\begin{proof}
	The proof can be found in \cite[Lemma 3.5]{GoHeMa22}. For the sake of completeness, it is also presented here.
	Let $\bar M$, $\bar L$ be the operators defined on the subgraph $\bar {\mathcal G}$ from \Cref{ass:network}. We consider the generalized eigenvalue problem $\bar L v = \lambda \bar M v$ posed in the space $\hat V(\bar {\mathcal N})$ with $\lambda_1 \leq \lambda_2 \leq \dots$ denoting its eigenvalues. Due to \Cref{ass:network}, the subgraph $\bar {\mathcal G}$ is connected and we have that $\lambda_1 = 0$ (corresponding to constant eigenfunctions) and $\lambda_2>0$. By the min-max theorem, the second eigenvalue admits the characterization
	\begin{equation*}
	\lambda_2 = \inf_{\substack{0\neq v \in V\\ \vspf{\bar M v}{1} = 0}} \frac{\vspf{\bar L v}{v}}{\vspf{\bar M v}{v}}.
	\end{equation*}
	Denoting with $c$ the $\bar M$-orthogonal projection of $v$ onto the constant functions, we obtain
	\begin{equation*}
	\mnormfo{v-c}{T}^2 \leq \mnormfo{v-c}{\mathsf N(T)}^2 \leq \lambda_2^{-1} \lnormfo{v-c}{\mathsf N(T)} ^2 =\lambda_2^{-1} \lnormfo{v}{\mathsf N(T)}^2
	\end{equation*} 
	which is the desired inequality.
\end{proof}

This lemma states, similarly as \Cref{l:friedrichs}, only the existence of a constant $C_\mathrm{po}$ and not its qualitative behavior in dependence of $H$. In \Cref{subsec:poincare}, a numerical experiment confirms that $C_\mathrm{po}$ is proportional to $H$ provided that the considered meshes are sufficiently coarse, see \Cref{ass:network}. Theoretically, such a scaling can be proved under the assumption of~a $d$-dimensional isoperimetric inequality, cf. \cite[Lemma~3.6]{GoHeMa22}. 

This yields us to the following assumption.

\begin{assumption}[Poincar\'{e} constant]
	\label{ass:poincareconstant}
	Assume that 
	there exists $\mu >0$ independent of $H$ such that 
	\begin{equation*}
	C_\mathrm{po} \leq \mu H.
	\end{equation*}
\end{assumption} 

\section{Prototypical approximation}\label{s:idealmethod}

In this section, we construct prototypical problem-adapted ansatz spaces with uniform approximation properties independent of the material data. The word prototypical shall emphasize that, without modification, the constructed ansatz spaces are not feasible in a practical implementation. 

A common technique for the construction of such problem-adapted ansatz spaces is the application of the inverse operator $\K^{-1}$ to some classical finite element spaces. Here, we adapt this technique to the spatial network setting. Let us consider the space of $\TH$-piecewise constants given by 
\begin{equation*}
\Pnull(\T_H)=\operatorname{span}\{\one_\elem\with \elem\in \mathcal T_H\} \subset \hat V,
\end{equation*}
with the indicator function equaling one for $x \in \mathcal N$ that are contained in $T$ and zero else. 

\subsection{$L^2$-projection}
We introduce the $L^2$-projection onto $\Pnull(\TH)$, defined as
\begin{equation}\label{eq:defL2}
\PiH\colon \hat V \to \Pnull(\TH),\quad v \mapsto \sum_{\elem \in \mathcal T_H}\frac{\vsp{M_T v}{1}}{\mnormfo{1}{T}^2}\one_\elem.
\end{equation}

The following lemma states global stability and approximation estimates for $\PiH$.  
\begin{lemma}[Properties of the $L^2$-projection]
	The $L^2$-projection is stable, i.e., for all $v \in \hat V$, it holds that
	\begin{equation}
	\label{eq:stabL2}
	\mnormf{\PiH v}\leq \mnormf{v}.
	\end{equation}
	Further, if \Cref{ass:network} and \Cref{ass:poincareconstant} are satisfied, then there is a constant $C_\Pi>0$ independent of $H$ such that the following  approximation estimate holds, for all~$v \in \hat V$,
	\begin{equation}
	\label{eq:appL2}
	\mnormf{v-\PiH v} \leq C_\Pi H \lnormf{v}.
	\end{equation}
\end{lemma}
\begin{proof}
	For the proof of \eqref{eq:stabL2}, we use the local stability of $\PiH$ which follow immediately from definition \eqref{eq:defL2}. This yields
	\begin{align*}
	\mnormf{\PiH v}^2 = \sum_{\elem \in \mathcal T_H} \mnormfo{(\PiH v)|_T}{T}^2 = \sum_{\elem \in \mathcal T_H} \vsp{M_T v}{1} \mnormfo{1}{T}^{-1} \leq \sum_{\elem \in \mathcal T_H}\mnormfo{v}{T}^2 = \mnormf{v}^2.
	\end{align*}
	For proving \eqref{eq:appL2}, we again split the norm into a sum of element contributions and employ \eqref{eq:poincare} and \Cref{ass:poincareconstant} to obtain
	\begin{align*}
	\mnormf{v-\PiH v}^2 = \sum_{\elem \in \mathcal T_H} \mnormfo{v-\PiH v}{T}^2\leq \mu^2 H^2 \sum_{\elem \in \mathcal T_H}\lnormfo{v}{\mathsf N(T)}^2 \leq 3^d\mu^2H^2\lnormf{v}^2,
	\end{align*}
	where the constant $3^d$ reflects the finite overlap of the patches $\mathsf N(T)$.
\end{proof}

\subsection{Prototypical method}
We define the prototypical problem-adapted ansatz space~as
\begin{equation}
\label{eq:idealsp}
V_H \coloneqq \K^{-1}\Pnull(\TH).
\end{equation}
The corresponding Galerkin method seeks a discrete approximation $u_H\in V_H$ such that, for all $v_H\in V_H$,
\begin{equation}\label{e:galerkinideal}
\vsp{K u_H}{v_H} = \vsp{M f}{v_H}.
\end{equation}

When using problem-adapted ansatz spaces as $V_H$, the approximation problem of the solution $u$ in $V_H$ is transformed into an approximation problem of the right-hand side $f$ in $\Pnull(\TH)$. This allows to show the desired uniform approximation properties without pre-asymptotic effects. 

\begin{lemma}[Uniform  approximation]\label{l:ua}
	Let the network satisfy \Cref{ass:network} and \Cref{ass:poincareconstant}. Then, for any $f \in \hat V$, the  prototypical approximation $u_H$ defined in \eqref{e:galerkinideal} converges quadratically in $H$, i.e., 
	\begin{equation}
	\label{eq:estidealmethod}
	\lnormf{u-u_H} 
	\leq C_\Pi \alpha^{-1}H \mnormf{f-\PiH f}
 \leq C_\Pi^2 \alpha^{-1}H^2\lnormf{f}.
	\end{equation}
\end{lemma}
\begin{proof}
	This is the spatial network counterpart of \cite[Lemma 3.2]{HaPe21b}. For the error estimate, we introduce $\bar u_H \coloneqq \K^{-1}\PiH f\in V_H$ and employ C\'ea's lemma for estimating the energy approximation error of $u_H$ against the one for $\bar u_H$, i.e., 
	\begin{equation*}
	 \knormf{u-u_H} \leq \knormf{u-\bar u_H}.
	\end{equation*}
	Denoting $e \coloneqq u-\bar u_H$, we further obtain using \eqref{eq:appL2}
	\begin{align*}
	 \knormf{e}^2  &= \vsp{Ke}{e} = \vsp{M(f-\PiH f)}{e} \leq \vsp{M(f-\PiH f)}{e - \PiH e}\\
	&\leq \mnormf{f-\PiH f}\mnormf{e - \PiH e}
	\leq C_\Pi\alpha^{-1/2} H\mnormf{f-\PiH f} |e|_K.
	\end{align*}
	Dividing by $\knormf{e}$, inferring that 
	\begin{equation*}
	\mnormf{f-\PiH f} \leq C_\Pi H \lnormf{f}
	\end{equation*}
	and using \eqref{eq:normeq}, the assertion follows.
\end{proof}

\begin{remark}[$\TH$-piecewise right-hand sides]
	\label{rem:exactnessideal}
	For $f \in \Pnull(\TH)$, the prototypical method is exact, as for $\TH$-piecewise constant right-hand sides, it holds $\mnormf{f-\PiH f} = 0$ in \eqref{eq:estidealmethod}.
\end{remark}

\section{Localization}\label{s:locapprox}

This section provides a localization strategy that turns the prototypical multiscale method introduced in the previous section into a practical method. 
Inspired by \cite{HaPe21b}, we introduce a localization strategy for spatial network models that identifies local $\TH$-piecewise constant source terms with rapidly decaying responses under the solution operator~$\K^{-1}$. This rapid decay enables a localization of the basis functions to element patches which paves the way to an efficiently computable localized ansatz space. 

In order to simplify the notation in the subsequent derivation, we fix an element $T \in \TH$ and its surrounding $\ell$-th order patch $\omega \coloneqq \mathsf N^\ell(T)$. Furthermore, let $\T_{H,\omega}$ denote the submesh of $\TH$ consisting of elements contained in $\omega$. 

\subsection{Localization ansatz}
For the construction of an almost local basis of the prototypical ansatz space~$V_H$ defined in \eqref{eq:idealsp}, we assign one basis function to each element $T \in \TH$. The (in general global) basis function $\psi \in V_H$ associated to element $T$ is determined by the following ansatz
\begin{equation}
\label{eq:psi}
\psi = \K^{-1} g\quad\text{ with  }\quad g \coloneqq \sum_{T \in \T_{H,\omega}} g_T\one_\elem \in \Pnull(\T_{H,\omega})
\end{equation}
with coefficients $(g_T)_{T\in \TH}$ to be determined subsequently. A local counterpart of $\psi$ which is supported on the patch $\omega$ can be defined by applying the patch-local solution operator~$\K_\omega^{-1}$ instead of  $\K^{-1}$, i.e., 
\begin{equation*}
\label{eq:phi}
\varphi \coloneqq \K_\omega^{-1} g.
\end{equation*}

In general, the locally supported function $\varphi$ is a poor approximation of its global counterpart $\psi$. However, we aim to choose the coefficients $(g_T)_{T\in \TH}$ such that $\psi$ is  well approximated by its patch-local counterpart $\varphi$. 

\subsection{Conormal derivatives}
For this purpose, we transfer the concept of conormal derivatives to the spatial network setting. Let $\tilde V_\omega$ denote the subset of $V$ consisting of functions that are supported on nodes in $\omega$ or its neighboring nodes. 

The following preliminary result is inevitable for the definition of conormal derivatives.
\begin{lemma}[Inner product on $\tilde V_\omega$]\label{l:lpmnorm}
	The bilinear form
	\begin{equation*}
	\vsp{(L_\omega + M_\omega)\, \cdot}{\cdot}
	\end{equation*}
	is an inner product on the space $\tilde V_\omega$.
\end{lemma}
\begin{proof}
	We only need to show the  positivity of the inner product which we do by contradiction. Assume that there exists $0 \neq v \in \tilde V_\omega$ such that
	$\vsp{(L_\omega +M_\omega) v}{v} = 0$. Let $\tilde{\mathcal G}$ denote the subgraph of $\mathcal G$ within $\omega$ extended by its neighboring nodes and the respective edges. We can decompose $\tilde{\mathcal G}$ into a finite number of connected components $\tilde{\mathcal G_{i}}$. For each connected component, we have that $v$ equals a constant $c_i$ (otherwise $(L_\omega v,v)>0$). Denoting with $M_i$ the mass matrix with respect to  $\tilde{\mathcal G_{i}}$, we have that $\vsp{M_i c_i}{c_i} = c_i^2\vspf{M_i 1}{1}$  which is only zero if $c_i = 0$. Here, it shall be noted that a connected component cannot have nodes solely outside of $\omega$ since, by definition, these nodes are connected to nodes within $\omega$. The assertion follows immediately.
\end{proof}

In the spatial network setting, we define the conormal derivative of $\varphi$, denoted by $B_\omega \varphi  \in \tilde V_\omega^\prime$, as the functional that  satisfies, for all $v \in \tilde V_\omega$,  
 \begin{equation}
 \label{eq:conormalder}
\vsp{B_{\omega} \varphi }{v}\ \coloneqq \vsp{K \varphi}{v} - \vsp{Mg}{v}.
\end{equation}
Further, we define its operator norm as
\begin{equation*}
|B_{  \omega} \varphi|_{\tilde V_{  \omega}^\prime} \coloneqq \sup_{v \in  \tilde V_{  \omega}} \frac{\vsp{B_{  \omega} \varphi}{v}}{\vnormo{v}{  \omega}}.
\end{equation*}
The following lemma shows that the operator norm of $B_\omega\varphi$ can be bounded in terms of~$\varphi$ and its corresponding right-hand side $g$.  
\begin{lemma}[Continuity of conormal derivative]
	The operator norm of $B_\omega \varphi$ can be bounded as follows 
	\begin{equation*}
	|B_\omega \varphi|_{\tilde V_\omega^\prime} \leq \max\{2\beta,1\}\sqrt{\lnormfo{\varphi}{\omega}^2 + \mnormfo{g}{\omega}^2}.
	\end{equation*}
\end{lemma}
\begin{proof}
	We denote by $\tilde {\mathcal G} = (\tilde {\mathcal N},\tilde {\mathcal E})$ the subgraph of $\mathcal G$ within $\omega$ extended by its neighboring nodes and the respective edges. We define a semi-norm on $\tilde V_\omega$ by
	\begin{equation*}
	|v|_{\tilde L,\omega}^2 \coloneqq \sum_{x \in \tilde {\mathcal N}} \frac{1}{2}\sum_{\tilde {\mathcal N} \ni y\sim x} \frac{(v(x)-v(y))^2}{|x-y|}
	\end{equation*}
	which is equivalent to $\lnormfo{\cdot}{\omega}$. More precisely, for all $v \in \tilde V_\omega$, it holds  $\lnormfo{v}{\omega} \leq |v|_{\tilde L,\omega}$ and
	\begin{align*}
	|v|_{\tilde L,\omega}^2 \leq 2 \sum_{x \in \mathcal N(\omega)} \frac{1}{2}\sum_{y\sim x} \frac{(v(x)-v(y))^2}{|x-y|} = 2 \lnormfo{v}{\omega}^2.
	\end{align*}
	Using that $g \in \hat V_\omega$, $\varphi \in V_\omega$, and $v \in \tilde V_\omega$, the norm equivalence, and the discrete Cauchy--Schwarz inequality, we obtain
	\begin{align*}
	\vsp{B_\omega \varphi }{v} &= \vspf{K \varphi }{v}-\vspf{M g}{\omega} \leq \beta |\varphi|_{\tilde L,\omega}|v|_{\tilde L,\omega} + \mnormfo{g}{\omega}\mnormfo{v}{\omega} \\
	&\leq 2\beta \lnormfo{\varphi}{\omega}\lnormfo{v}{\omega} + \mnormfo{g}{\omega}\mnormfo{v}{\omega}\\
	&\leq \max\{2\beta,1\}\sqrt{\lnormfo{\varphi}{\omega}^2 + \mnormfo{g}{\omega}^2}\vnormfo{v}{\omega}
	\end{align*}
	which is the boundedness of the conormal derivative.
	\end{proof}

The following result enables a computation of the $\tilde V_{ \omega}^\prime$-norm of the conormal derivative and is key for the method's practical implementation.
\begin{lemma}[Computation of the conormal derivative's norm]
	\label{l:conormalder}
	The operator norm of the conormal derivative $B_{  \omega} \varphi \in \tilde V_{ \omega}^\prime$ can be computed as
\begin{equation*}
|B_{  \omega} \varphi|_{\tilde V_{  \omega}^\prime} = \vnormfo{\tau}{\omega},
\end{equation*}
where $\tau \in \tilde V_\omega$ solves, for all $v \in \tilde V_{  \omega}$,
	\begin{equation}
	\label{eq:tau}
	 \vsp{(L_{  \omega} + M_\omega) \tau }{v} = \vsp{B_{  \omega} \varphi }{v}.
	\end{equation}
\end{lemma}
\begin{proof}
	\Cref{l:lpmnorm} yields the unique existence of a solution $\tau$ to problem \eqref{eq:tau}. 
	Furthermore, using~\eqref{eq:tau}, we can compute the operator norm of $B_\omega \varphi$ as follows
	\begin{align*}
		|B_{  \omega} \varphi|_{\tilde V_{  \omega}^\prime} \coloneqq \sup_{v \in \tilde V_{  \omega}} \frac{\vsp{B_{  \omega} \varphi}{v}}{\vnormo{v}{  \omega}} = \sup_{v \in \tilde V_{ \omega}} \frac{\vsp{L_{  \omega} \tau }{v} + \vsp{M_{  \omega} \tau}{v}}{\vnormfo{v}{  \omega}} = \vnormfo{\tau}{  \omega},
	\end{align*}
	where we used that the supremum is attained for $v = \tau$. 
\end{proof}

\subsection{Choice of local basis}
It turns out that  $\varphi$ approximates $\psi$ well  provided that the parameters $(g_T)_{T \in \TH}$ are chosen such that the conormal derivative of $\varphi$  is small, since, for all $v \in V$, it holds
\begin{equation*}
\label{eq:keyobservation}
\vsp{K(\varphi-\psi)}{v} = \vsp{K \varphi}{v}-\vsp{M g}{v} = \vsp{B_{ \omega} \varphi}{v}.
\end{equation*}

\begin{figure}[h]
	\includegraphics[width=.3\linewidth]{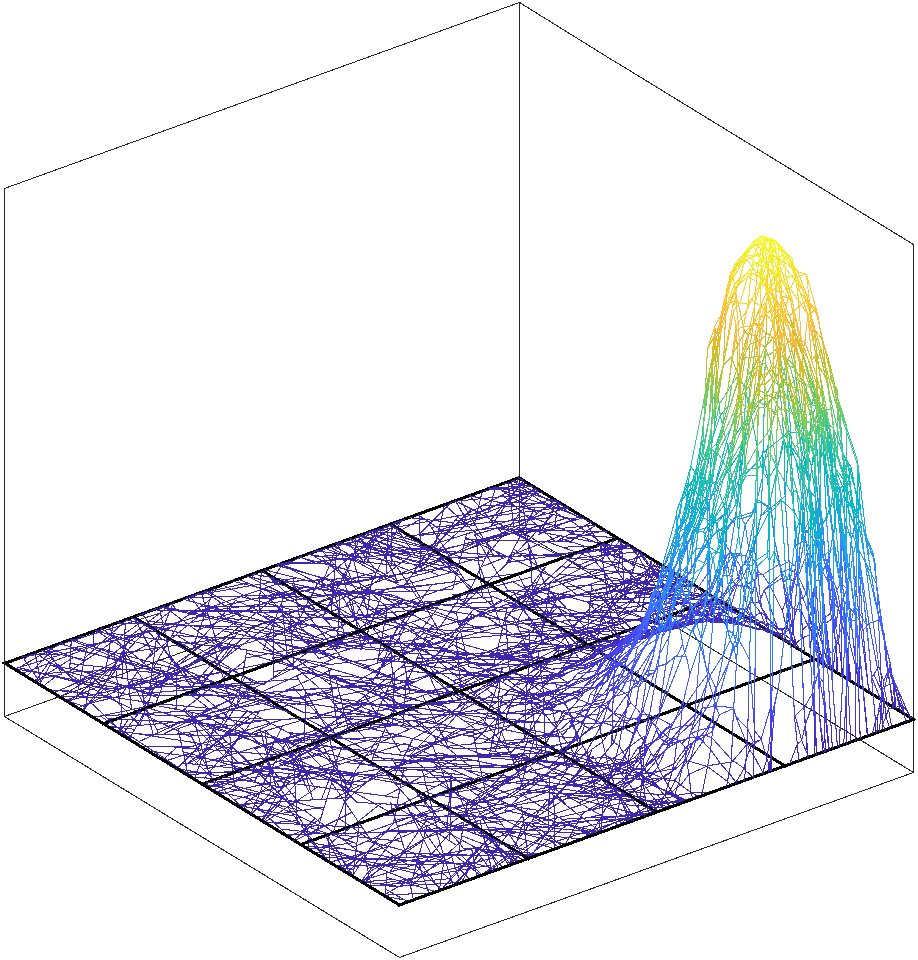}\hfill
	\includegraphics[width=.3\linewidth]{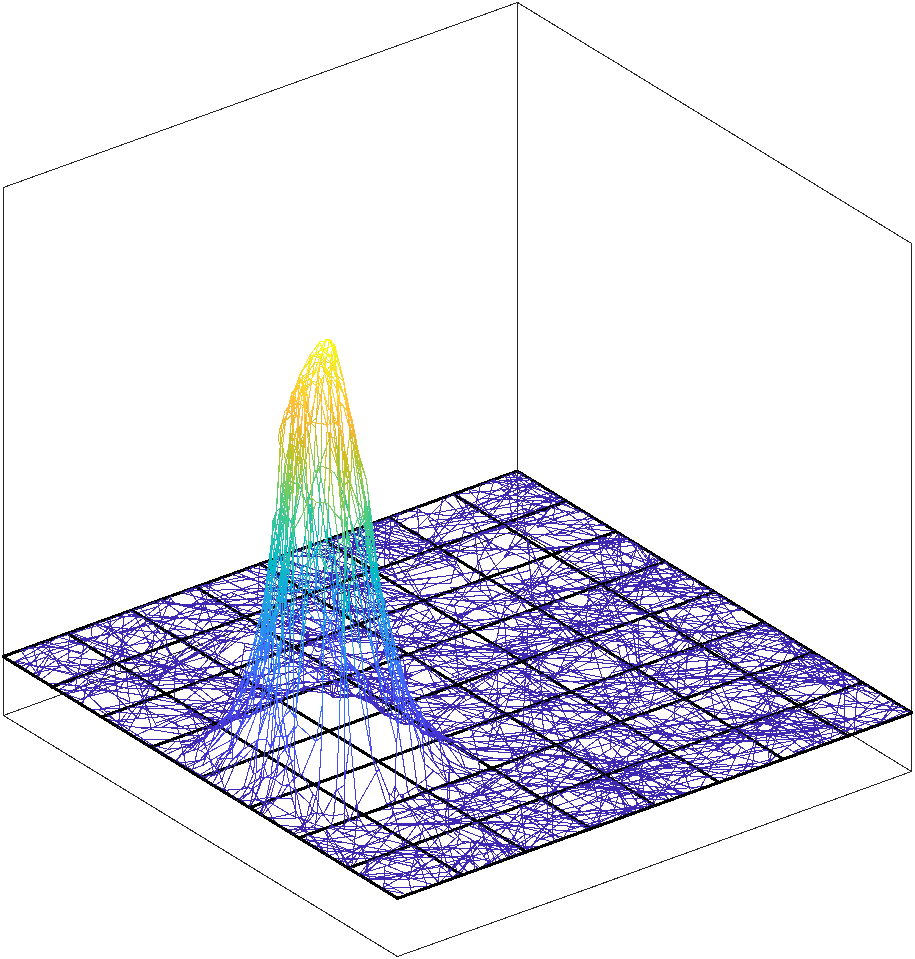}\hfill
	\includegraphics[width=.3\linewidth]{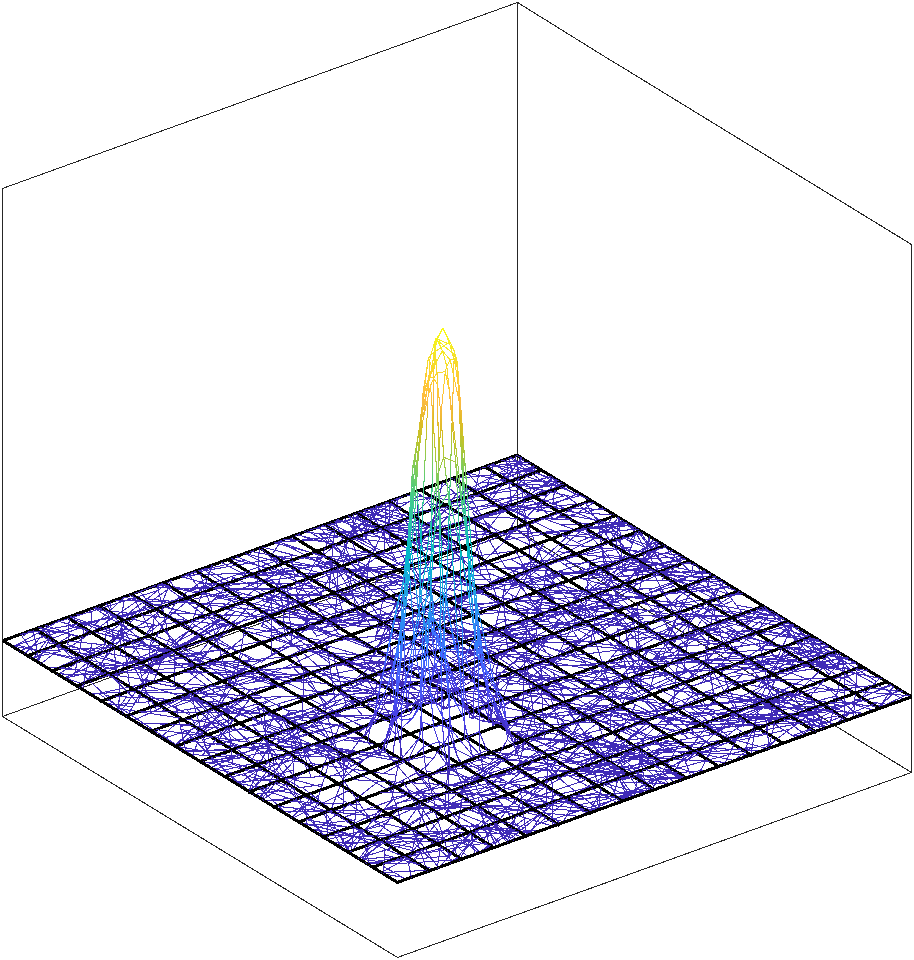}\vspace{.5cm}
	\includegraphics[width=.3\linewidth]{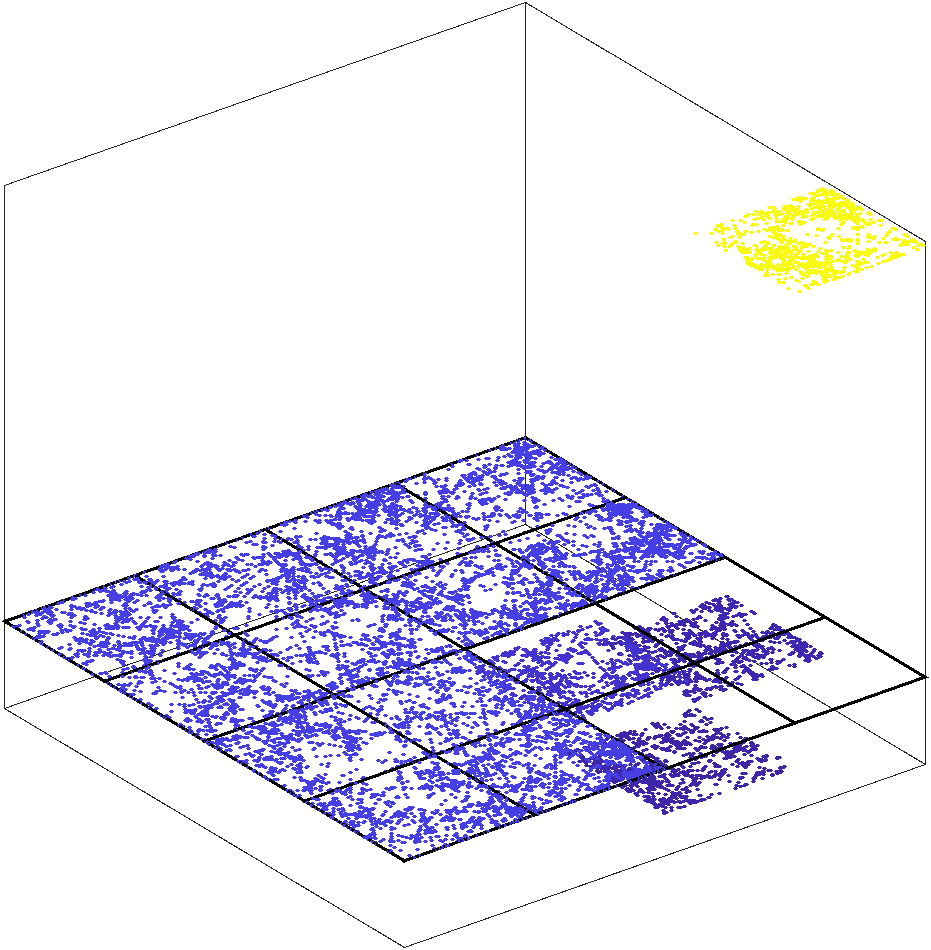}\hfill
	\includegraphics[width=.3\linewidth]{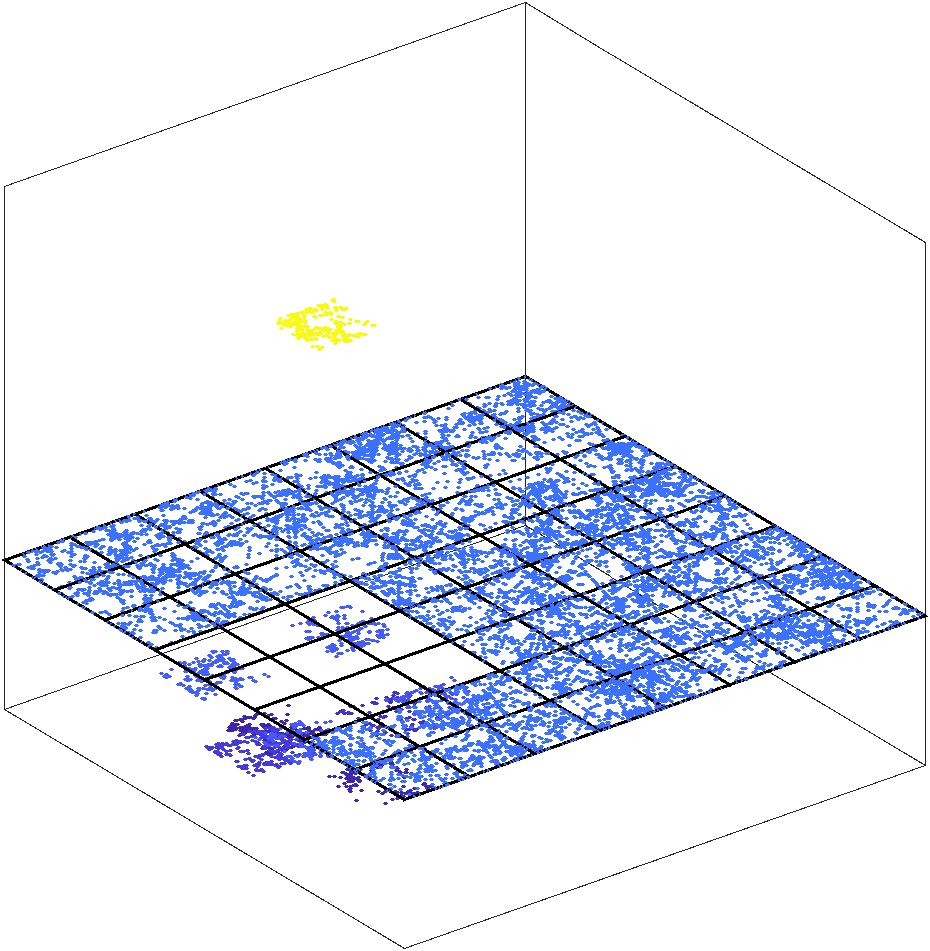}\hfill
	\includegraphics[width=.3\linewidth]{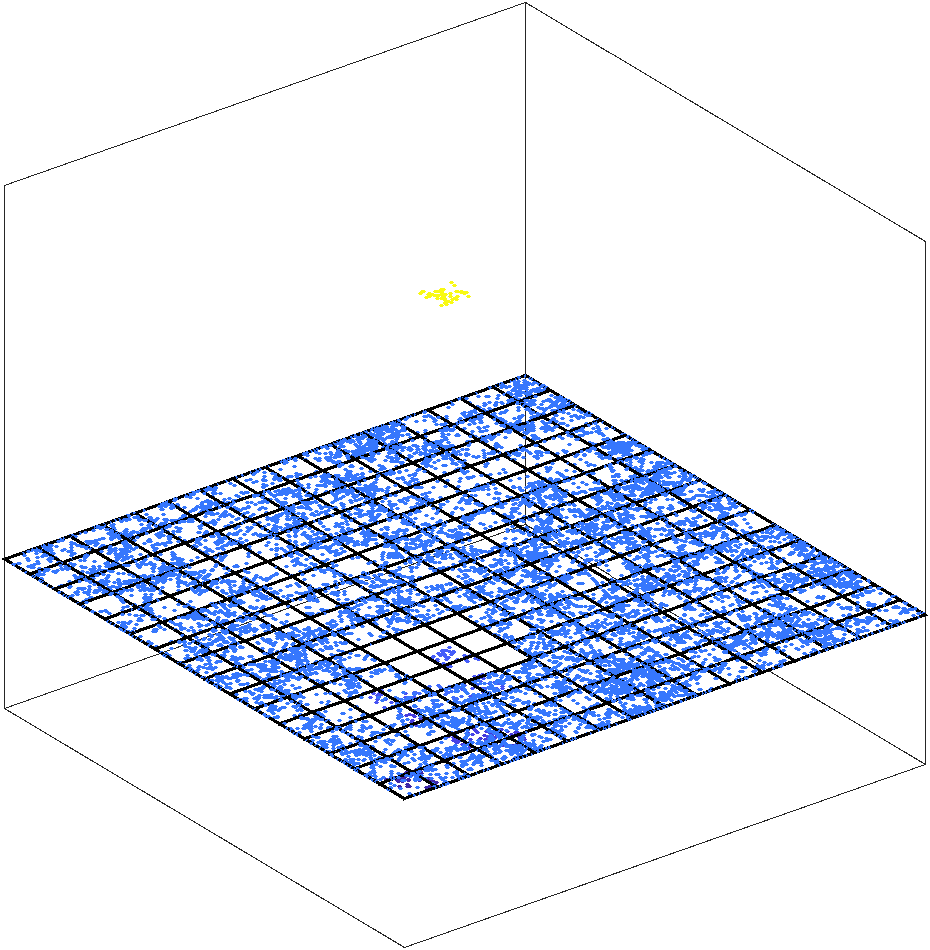}
	\caption{Basis functions $\varphi$ for $\ell = 1$ and several coarse mesh sizes $H$ (top) and corresponding right-hand sides $g$ (bottom).}
	\label{fig:basisfuns}
\end{figure}

Hence, denoting by $\R \colon \Pnull(\T_{H,\omega}) \to \tilde V_\omega$, $g\mapsto \tau$ the linear mapping that maps the right-hand side $g$ to its corresponding function $\tau$ defined in \eqref{eq:tau}, we choose $g$ as an  $\mnormfo{\cdot}{\omega}$-normalized minimizer of the following quadratic constraint minimization problem 
\begin{equation}
\label{eq:g}
g \in \argmin_{q \in \Pnull(\T_{H,\omega})}\frac{\vnormfo{\R q}{\omega}^2}{\mnormfo{q}{\omega}^2}.
\end{equation}
Due to \Cref{l:conormalder}, this is equivalent to the minimization of the conormal derivative. For a depiction of selected basis functions and their corresponding right-hand sides, we refer to \Cref{fig:basisfuns}.

\begin{remark}[Numerical implementation]
	Numerically, instead of solving \eqref{eq:g}, one can equivalently solve for the eigenvector corresponding to the smallest eigenvalue of the low-dimensional generalized eigenvalue problem 
	\begin{equation}
	\label{eq:evp}
	Ax = \lambda Cx
	\end{equation}
	with matrices $A,C \in \mathbb R^{N\times N}$, $N \coloneqq \# \T_{H,\omega}$, defined as
	\begin{align*}
	A_{ij} = \vspf{(L_\omega + M_\omega) \R \one_{\elem_j}}{\R \one_{\elem_i}},\quad C_{ij} \coloneqq \vspf{M_\omega \one_{\elem_j}}{\one_{\elem_i}},
	\end{align*}
	where $\{T_i\with i = 1,\dots N\}$ is some numbering of the elements in $\T_{H,\omega}$.
\end{remark}

As the notation in \eqref{eq:g} already indicates, the solution to the minimization problem is, in general, non-unique. Indeed, for large oversampling parameters $\ell$ and for patches~$\omega$ close to the boundary $\partial \Omega$, there might by several (linearly independent) $g$ with similar in size Rayleigh quotients. In such cases, the appropriate choice of $g$ can be difficult, cf. \Cref{subsec:stabrhs} for an implementation resolving this issue in practice.

\subsection{Localization error}
We define the quantity $\sigma_T$ as the norm of the conormal derivative of the selected basis function
\begin{equation}
\label{eq:sigma}
\sigma_T = \sigma_T(H,\ell) = \vnormfo{\R g}{\omega}.
\end{equation}
The quantity $\sigma_T$ determines the local localization error, i.e., the error when approximating $\psi$ by its local counterpart $\varphi$.

\section{Super-localized approximation}\label{s:locmethod}

Using the localized basis function introduced in the previous section, this section transforms the prototypical method \eqref{e:galerkinideal} into a practically feasible method. For a fixed oversampling parameter, we define the localized ansatz space as
\begin{equation*}
V_{H,\ell} \coloneqq \linh\{\varphi_{\elem,\ell}\with \elem \in \mathcal T_H\}\subset V.
\end{equation*}
Note that, in the case of ambiguity, we write $\varphi_{\elem,\ell}$, $\psi_{\elem,\ell}$, and $g_{\elem,\ell}$ for $\varphi$, $\psi$, and $g$ in order to emphasize their dependence on $T,\ell$.

The SLOD determines the Galerkin approximation of the solution $u$ to \eqref{eq:weakform} in the space $V_{H,\ell}$, i.e., it  seeks $u_{H,\ell}\in V_{H,\ell}$ such that, for all $v_{H,\ell}\in V_{H,\ell}$,
\begin{equation}\label{e:wflocmethod}
\vspf{K u_{H,\ell}}{v_{H,\ell}} = \vspf{M f}{v_{H,\ell}}.
\end{equation}

\subsection{Riesz basis property of right-hand sides}\label{subsec:stabrhs}
The choice of right-hand sides \eqref{eq:g}, in general, does not guarantee their linear independence. For the analysis, we assume that the right-hand sides $\{g_{\elem,\ell}\with \elem \in \TH\}$ span $\Pnull (\mathcal T_H)$ in a stable way. Subsequently, we also present an implementation strategy that ensures the basis stability in practice.

\begin{assumption}[Riesz stability]\label{a:Rieszbasis}
The set $\{g_{\elem,\ell}\with \elem \in \TH\}$ is a Riesz basis of  $\Pnull(\mathcal T_H)$, i.e., there is $C_{\mathrm{r}}(H,\ell)>0$ depending polynomially on $H, \ell$ such that, for all $(c_\elem)_{T \in \TH}$,
	\begin{equation}
	\label{eq:riesz}
	C^{-1}_{\mathrm{r}}(H,\ell)\sum_{\elem \in \mathcal T_H} c_\elem^2  \leq \Big\vert \sum_{\elem \in \mathcal T_H} c_\elem g_{\elem,\ell}\Bigl\vert_{M}^2 \;\leq C_{\mathrm{r}}(H,\ell) \sum_{\elem \in \mathcal T_H} c_\elem^2.
	\end{equation}
\end{assumption}

The Riesz basis property is closely related to the eigenvalues of the matrix containing inner products of the right-hand sides $g_{T,\ell}$ as the following remark shows.
 \begin{remark}[Riesz constant]
 	The constant $C_\mathrm{r}$ is determined by the smallest and largest eigenvalue of the matrix $G \in \mathbb R^{M\times M}$, $M\coloneqq \#\TH$ defined as
 	\begin{equation*}
 	G_{ij} = \vspf{Mg_{T_j,\ell}}{g_{T_i,\ell}},
 	\end{equation*}
 	where $\{T_i\with i = 1,\dots,M\}$ is some numbering of the elements in $\TH$. Denoting its eigenvalues by $\lambda_1\leq \lambda_2 \leq \dots, \leq \lambda_M$, the constant in the lower (resp. upper) bound  of~\eqref{eq:riesz} is then the smallest (resp. largest) eigenvalue of $G$. Thus, we can set
 	$$C_\mathrm{r} = \max\{\lambda_M,\lambda_1^{-1}\}.$$
 \end{remark}

The stability issue of the basis of right-hand sides similarly appears in the continuous setting for the SLOD from \cite{HaPe21b}. Therein, an algorithm is proposed that selects the right-hand sides in a stable and efficient manner.
The algorithm is based on the observation that stability issues only occur for patches close to the boundary $\partial \Omega$. By grouping certain troubled patches and solving one minimization problem of the  type \eqref{eq:g} for such a group of patches, the stability of the right-hand sides within this group of patches can be ensured by construction. As this procedure only needs to be applied close to the global boundary for certain groups of patches, the algorithm only introduces minimal communication between the patches. For more information and a detailed illustrative description of the algorithm, see \cite[Appendix B]{HaPe21b}.

\subsection{A-posteriori error estimate}
We derive an error estimate for the SLOD solution~\eqref{e:wflocmethod} which is explicit in the quantity 
\begin{equation}\label{e:qom}
\sigma = \sigma(H,\ell) \coloneqq \max_{\elem\in\TH} \sigma_\elem(H,\ell)
\end{equation}
with $\sigma_T$ defined in  \eqref{eq:sigma}. The quantity $\sigma$ determines the size of the method's localization error.  \Cref{rem:decaysigma} presents a summary of existing decay results for $\sigma$ in the continuous setting and the spatial network setting.

\begin{theorem}[Uniform localized approximation]\label{t:error}
		Let the network satisfy \Cref{ass:network} and \Cref{ass:poincareconstant}. Further, let $\{g_{\elem,\ell}\with \elem \in \TH\}$ be stable in the sense of \Cref{a:Rieszbasis}. Then, for any $f \in \hat V$, the SLOD approximation defined in \eqref{e:wflocmethod} converges quadratically in~$H$ plus a localization error, i.e., there exists $C,C^\prime>0$ such that, for all $f \in \hat V$,
	\begin{equation}
	\label{eq:est}
	\begin{aligned}
	\lnormf{u-{u_{H,\ell}}}&\leq   C\hphantom{^\prime}\big(H \mnormf{f-\Pi_H f} + C_{\mathrm{r}}^{1/2}(H,\ell)\ell^{d/2} \sigma(H,\ell)\mnormf{f}\big)\\
	 &\leq C^\prime\big(H^2 \lnormf{f} + C_{\mathrm{r}}^{1/2}(H,\ell)\ell^{d/2} \sigma(H,\ell)\mnormf{f}\big)
	\end{aligned}
	\end{equation}
	with $C_{\mathrm{r}}$ from \Cref{a:Rieszbasis}. 
\end{theorem}

\begin{proof}
	Using C\'ea's Lemma, we can estimate the energy error of $u_{H,\ell}$ approximating $u$ by the respective energy error for $v_{H,\ell}$, for any $v_{H,\ell} \in V_{H,\ell}$. This yields
	\begin{align*}
	\knormf{u-u_{H,\ell}}\leq \knormf{u-v_{H,\ell}}.
	\end{align*}
	Next, we add and subtract $\bar u_H \coloneqq  \mathcal K^{-1}\PiH f$ and obtain with the triangle inequality
	\begin{align*}
	\knormf{u-u_{H,\ell}}\leq  \knormf{u-\bar u_H}+\knormf{\bar{u}_H-v_{H,\ell}}.
	\end{align*}
	
	The first term has already been estimated in the proof of \Cref{l:ua}. Therein, it was shown that
	\begin{equation*}
		\knormf{u-\bar u_H} 
	\leq C_\Pi^2 \alpha^{-1/2}H^2\lnormf{f}.
	\end{equation*}
	
	For the second term, we use that the prototypical method \eqref{e:galerkinideal} is exact for piecewise constant right-hand sides, i.e., in particular for $\PiH f$. This enables us to represent $\bar u_H$ using the functions $\psi_{T,\ell} = \K^{-1} g_{T,\ell}$ from \eqref{eq:psi} as
	\begin{equation*}\label{e:u_h}
	\bar u_H = \sum_{\elem \in \mathcal T_H} c_\elem\, \psi_{\elem,\ell},
	\end{equation*}
	where $(c_\elem)_{T \in \TH}$ are the coefficients of the expansion of $\PiH f$ in terms of the right-hand sides~$g_{\elem,\ell}$. This is possible as the $g_{T,\ell}$ are linearly independent which is, in particular, guaranteed by \Cref{a:Rieszbasis}. 
	For the particular choice \begin{equation*}
	v_{H,\ell} \coloneqq \sum_{\elem\in \mathcal T_H} c_\elem\, \varphi_{\elem,\ell}\in V_{H,\ell},
	\end{equation*}
	we obtain defining $e\coloneqq\bar u_H - v_{H,\ell}$ 
	\begin{align}
	\label{e:proof1}
	\knormf{e}^2 = \sum_{\elem \in \mathcal T_H}c_\elem\, \vspf{K(\psi_{\elem,\ell}-\varphi_{\elem,\ell})}{e}.
	\end{align}
	Writing $\mathsf N^\ell(T)$ instead of $\omega$, we obtain using \Cref{l:conormalder} and the definitions of $\psi_{\elem,\ell}$, the conormal derivative $B_{\mathsf N^\ell(T)}$ in \eqref{eq:conormalder}, and $\sigma_T$ from \eqref{eq:sigma} that
	\begin{equation}\label{e:proof2}
	\begin{aligned}
	\vspf{K(\psi_{\elem,\ell}-\varphi_{\elem,\ell})}{e}  &= \vspf{M g_{T,\ell}}{e} - \vspf{K \varphi_{\elem,\ell}}{e}\\& = -\vspf{B_{\mathsf N^\ell(T)} \varphi_{\elem,\ell}}{e} = - \vspf{B_{\mathsf N^\ell(T)}\varphi_{\elem,\ell}}{\tilde e} \\&\leq \sigma_\elem(H,\ell)\vnormfo{\tilde e}{\mathsf N^\ell(T)}  = \sigma_\elem(H,\ell)\vnormfo{e}{\mathsf N^\ell(T)}.
	\end{aligned}
	\end{equation}
	Here, we extended definition \eqref{eq:conormalder} to test functions in $V$ and use the notation  $\tilde{e}$ for the restriction of $e \in V$ to the space $\tilde V_{\mathsf N^\ell(T)}$. In the last step, we utilized $\vnormfo{\tilde e\,}{\mathsf N^\ell(T)} = \vnormfo{e}{\mathsf N^\ell(T)}$ which holds as the norm only considers nodes in the support of $\tilde e$.
	
	The combination of \eqref{e:proof1}, \eqref{e:proof2}, \Cref{a:Rieszbasis}, the discrete Cauchy--Schwarz inequality, the finite overlap of the patches, and \eqref{eq:stabL2} yields 
	\begin{equation*}
	\begin{aligned}
	\knormf{e}^2 &= \sum_{\elem \in \mathcal T_H}c_\elem\, \vspf{K(\psi_{\elem,\ell}-\varphi_{\elem,\ell})}{e}\leq  \sigma(H,\ell)\sum_{\elem \in \mathcal T_H}c_\elem\vnormfo{e}{\mathsf N^\ell(T)} \\& \leq \sigma(H,\ell)\sqrt{\sum_{\elem \in \mathcal T_H}c_\elem^2}\sqrt{\sum_{\elem \in \mathcal T_H} \vnormfo{e}{\mathsf N^\ell(T)}^2}\\
	&\leq  \sigma(H,\ell)\,C_{\mathrm{r}}^{1/2}(H,\ell)\mnormf{\PiH f}\,C_\mathrm{ol}\ell^{d/2}\vnormf{e} \\
	&\leq C_\mathrm{ol} C_{\mathrm{r}}^{1/2}(H,\ell)\ell^{d/2}\sigma(H,\ell)\alpha^{-1/2}\sqrt{1+C_\mathrm{fr}^2}\knormf{e} \mnormf{f}
	\end{aligned}
	\end{equation*}
	with constant $C_\mathrm{ol}>0$ reflecting the overlap of the patches $\mathsf N^\ell(T)$. In the last step, we applied Friedrichs' inequality \eqref{eq:friedrichs} on the full network. Putting together the previous estimates and using \eqref{eq:edgeweights}, the assertion can be concluded.  
\end{proof}

\begin{remark}[Decay of localization error]\label{rem:decaysigma}
	In the continuous setting, \cite{HaPe21b} conjectures that $\sigma$ decays super-exponentially in $\ell$. More precisely, there exists $C_\sigma>0$ depending polynomially on $H,\ell$ and $C_\mathrm{d}>0$ independent of $H,\ell$ such that, for all $\ell$,
	\begin{equation}
		\label{eq:sigmadec}
	\sigma(H,\ell) \leq C_\sigma(H,\ell) \exp(-C_\mathrm{d} \ell^\frac{d}{d-1}).
	\end{equation}
	In \cite[Section~7]{HaPe21b}, this result is justified theoretically using a conjecture from spectral geometry. Numerical experiments in \cite[Section~8]{HaPe21b} confirm the super-exponential decay numerically. Furthermore,  \cite[Lemma~6.4]{HaPe21b} provides a rigorous proof that $\sigma$ decays at least exponentially in $\ell$.
	
	In the spatial network setting, we also observe a rapid decay of the quantity $\sigma$ as $\ell$ is increased. Qualitatively, the decay behavior is is similar to the one in the continuous setting, see the numerical experiments in  \Cref{subsec:sigma} and \Cref{subsec:loc}. Using techniques from \cite{EGHKM22}, one can show, similarly as in the continuous setting, that $\sigma$ decays at least exponentially.
\end{remark}

\begin{remark}[$\TH$-piecewise right-hand sides]
	\label{rem:exactness}
	For $f \in \Pnull(\TH)$, only the localization error in \eqref{eq:est} is present, as for $\TH$-piecewise constant right-hand sides, it holds  $\mnormf{f-\PiH f} = 0$.
\end{remark}
\section{Numerical experiments}\label{s:ne}
For the subsequent numerical experiments, we utilize a spatial network constructed as follows. First, we sample 20,000 lines of length 0.05 which are uniformly rotated and with midpoints uniformly distributed in $[-0.025,1.025]^2$. Next, we remove all line segments outside of the unit square so that all lines are contained in the domain $\Omega = [0,1]^2$. We then define the network nodes as the line segments' endpoints and intersections. Two nodes are connected by an edge if the nodes share a line segment. We only consider the largest connected component of the network in order to ensure that the network is connected. We also remove all hanging nodes (nodes of degree one) in the interior of the domain along with the respective edges. All nodes at the boundary $\partial \Omega$ are Dirichlet nodes. The total number of nodes is around $80,000$.

\subsection{Poincar\'{e} constant}\label{subsec:poincare}

This numerical experiment shall justify \Cref{ass:poincareconstant} numerically. Given $T \in  \TH$, $H \in \H$, we construct the subgraph $\bar{\mathcal G} = (\bar {\mathcal N},\bar {\mathcal E})$ by means of a breadth-first search algorithm. Then, we compute the second eigenvalue $\lambda_2$ of the generalized eigenvalue problem $\bar L v = \lambda \bar M v$ posed in the space $\hat V(\bar{\mathcal N})$ with $\bar L$, $\bar M$ being defined on $\bar {\mathcal G}$,  see also the proof of \Cref{l:poincare}.

\begin{figure}[h]
	\includegraphics[width=.49\linewidth]{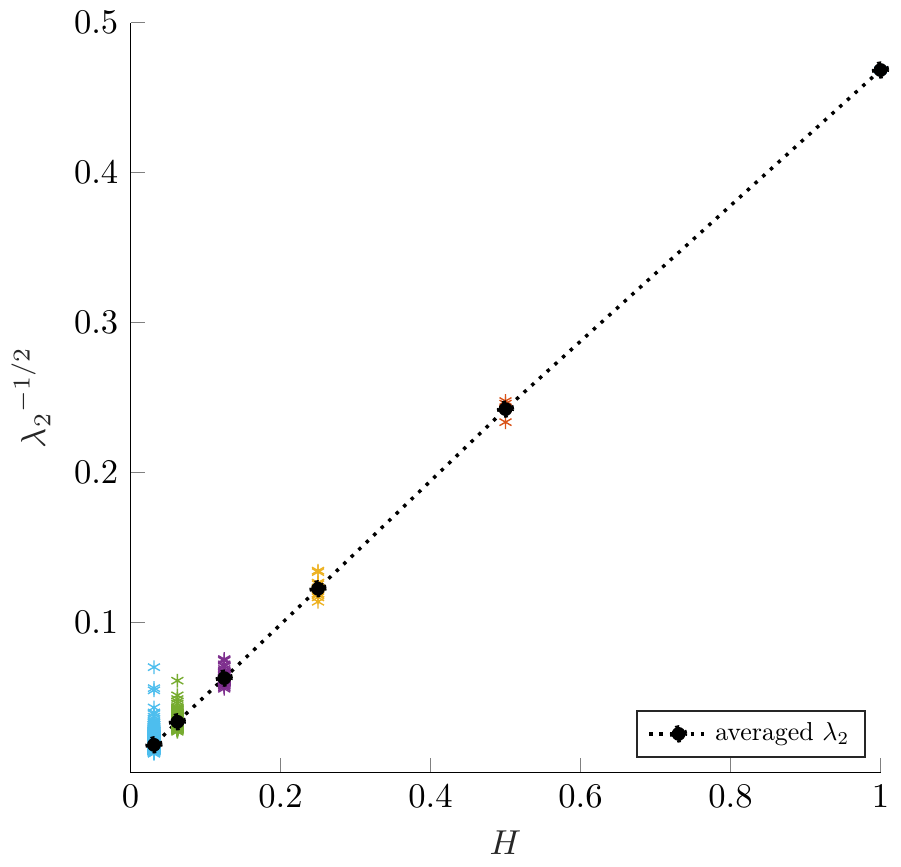}
	\caption{Eigenvalues $\lambda_2$ for different mesh sizes $H$.}
	\label{fig:eigenvalues}
\end{figure}

In \Cref{fig:eigenvalues}, the reciprocal square root of the numerically computed eigenvalues $\lambda_2$ corresponding to elements $T \in \TH$ is depicted for different mesh sizes $H \in \H$. The dotted black line interpolates, for all mesh sizes $H$, the averaged eigenvalue for the respective $H$. Note that the scaling of the $y$-axis is chosen such that the black line should appear linear if the desired scaling from \Cref{ass:poincareconstant} holds. 

\Cref{fig:eigenvalues} confirms \Cref{ass:poincareconstant} numerically. We note that for smaller mesh sizes, there is a considerable variation in the eigenvalues. This happens when we reach the critical mesh size, which is $H_0\approx 2^{-5}$ in this example.  

\subsection{Decay of $\sigma$}\label{subsec:sigma}
\begin{figure}[h]
	\includegraphics[width=.49\linewidth]{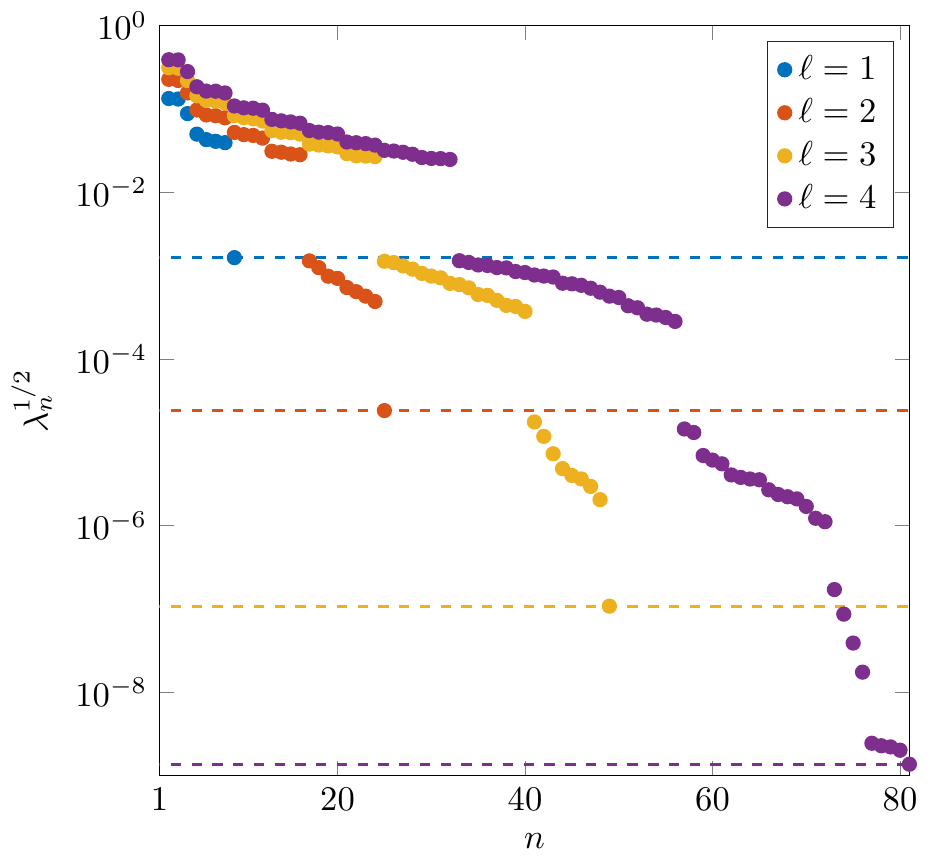}
	\caption{Eigenvalues of the eigenvalue problem \eqref{eq:evp} for patches of orders $\ell$. The dashed lines marks the values of $\sigma_T$ for the respective patches.}
	\label{fig:sigma}
\end{figure}

In this experiment, we numerically investigate the decay of $\sigma$ defined in~\eqref{e:qom} which is the maximum of the local quantities $\sigma_\elem$ from \eqref{eq:sigma}. For illustration purposes, we pick an element $T \in \T_{2^{-4}}$ whose fourth order patch has no intersection with the global boundary $\partial \Omega$.
\Cref{fig:sigma} then depicts the square root of the eigenvalues of the eigenvalue problem \eqref{eq:evp} for the patches $\mathsf N^\ell(T)$ with $\ell = 1,\dots,4$. By definition, the square root of the smallest eigenvalue coincides with $\sigma_T$.  The values of $\sigma_T$ are marked using dashed horizontal lines.

In \Cref{fig:sigma}, one observes a rapid decay of $\sigma_T$ as $\ell$ is increased. Note that for $\ell \geq 3$, the difference in magnitude of the smallest and largest eigenvalue of \eqref{eq:evp} is at least of order~$10^{14}$ which means that the respective matrices have a large condition number. This affects the accuracy of the eigenvalue solver and explain the flatting of the decay for large oversampling parameters as observed in \Cref{fig:sigma} and \Cref{fig1}.

\subsection{Super-localization}\label{subsec:loc}

\begin{figure}[h]
	\includegraphics[width=.49\linewidth]{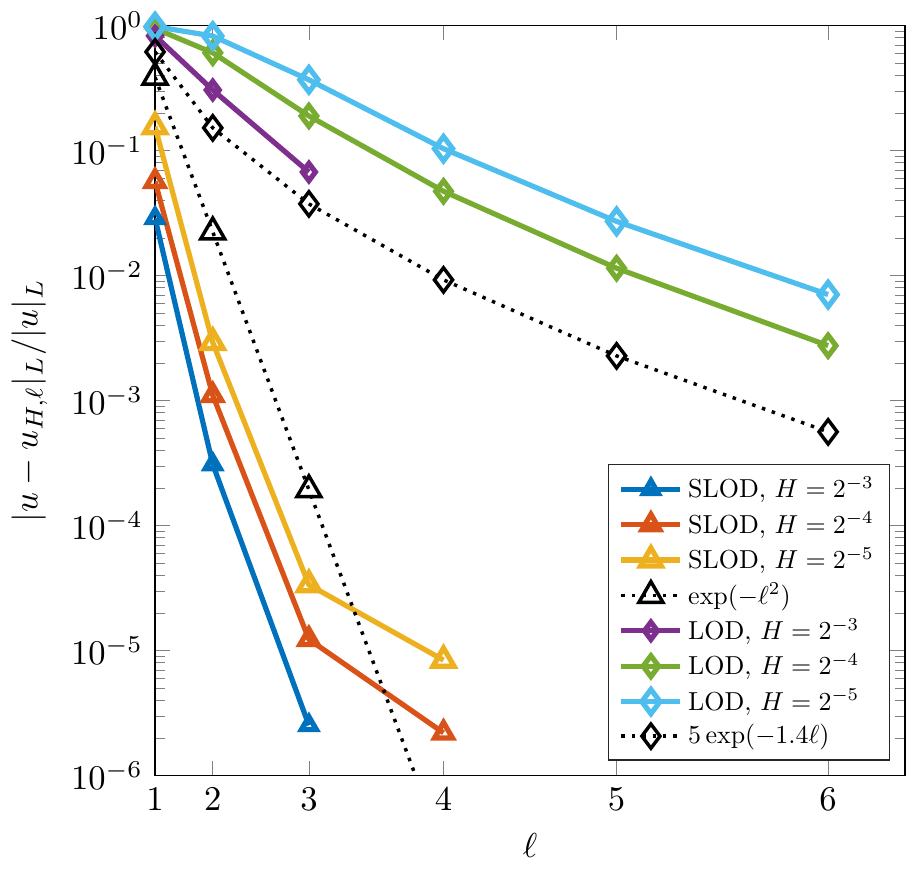}\hfill
	\includegraphics[width=.49\linewidth]{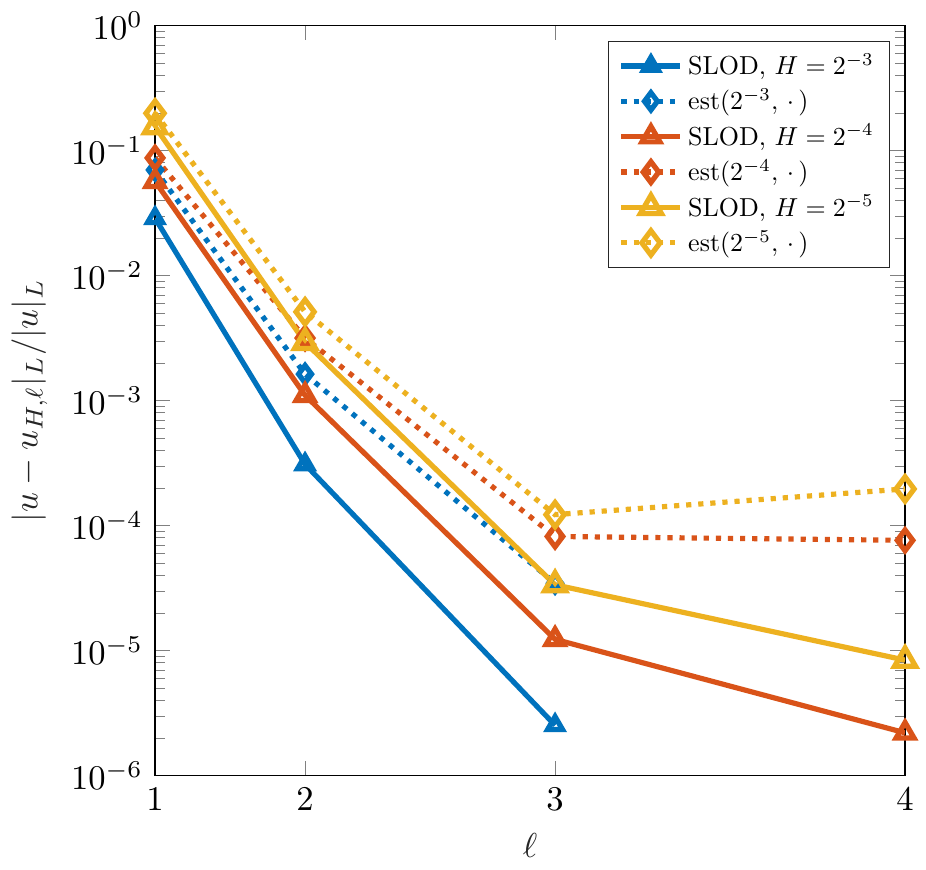}
	\caption{Plot of the relative localization errors in dependence of the oversampling parameter $\ell$ of the SLOD and the LOD (left). Depiction of the relative localization errors of the SLOD with the estimator \eqref{eq:estimator} (right).}
	\label{fig1}
\end{figure}

For this numerical experiment, we consider a weighted graph Laplacian  \eqref{eq:wgraphlaplacian} with edge weights $\gamma_{xy}$ independently and uniformly distributed in the interval $[0.01,1]$. We choose the right-hand side $f \equiv 1$, as for $f\in\Pnull(\TH)$, the error is bounded solely by the localization error, see \Cref{rem:exactness}. In \Cref{fig1} (left), the localization errors for the SLOD are shown. The localization errors are plotted for several coarse grids $\TH$ with respect to $\ell$. We additionally depict the localization errors when using the LOD for localizing the same prototypical ansatz space $V_H$ from \eqref{eq:idealsp}. As reference, we indicate lines showing the expected rates of decay of the localization errors for the SLOD and LOD which is super-exponential decay for the SLOD, cf. \eqref{eq:sigmadec}, and exponential decay for the LOD.
In \Cref{fig1} (right), we again depict the localization errors for the SLOD but this time together with the values of its estimator
\begin{equation}
\label{eq:estimator}
\mathrm{est}(H,\ell) \coloneqq C_\mathrm{r}^{1/2}(H,\ell)\ell^{d/2}\sigma(H,\ell)
\end{equation}
appearing in the error estimate in \Cref{t:error}. Note that the scaling of the $x$-axis in \Cref{fig1} is chosen such that a super-exponentially decaying curve appears to be a straight~line. 

\Cref{fig1} (left) numerically confirms the super-exponential decay rates of the localization errors as known for the continuous setting \cite{HaPe21b}. Note that, similarly as for the decay of $\sigma$ in \Cref{subsec:sigma}, we can also observe a flattening of the decay for $\ell \geq 3$. This might again be explained by the high condition number of the matrices in \eqref{eq:evp}, see \Cref{subsec:sigma}. 
The localization error of the LOD, depicted in \Cref{fig1} (left), decays exponentially, see e.g. \cite{Peterseim2021,EGHKM22}. Much larger values of $\ell$ are necessary in order to reach the accuracy level of the SLOD. 

\Cref{fig1} (right) shows that the error estimator is quite well aligned with the localization error and thereby underlines the sharpness of the estimator. For $\ell = 4$, we observe that also the estimator is slightly affected by the aforementioned conditioning issue.

\subsection{Optimal convergence}
\begin{figure}[h]
	\includegraphics[width=.499\linewidth]{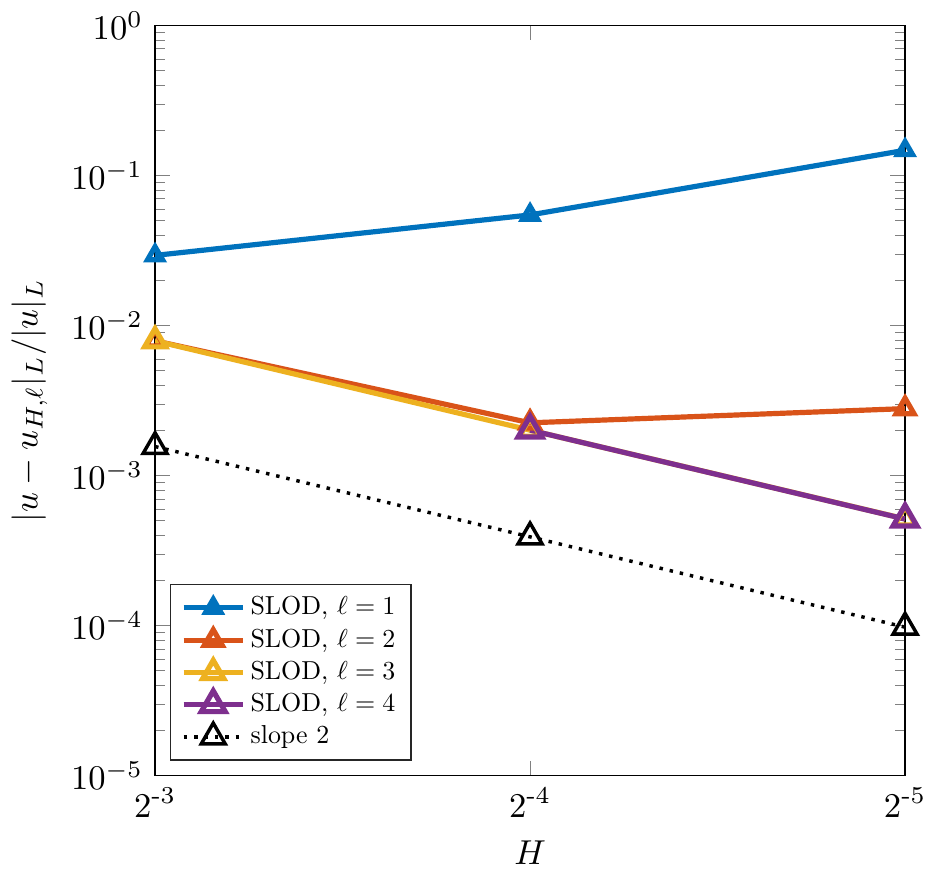}\hfill
	\includegraphics[width=.499\linewidth]{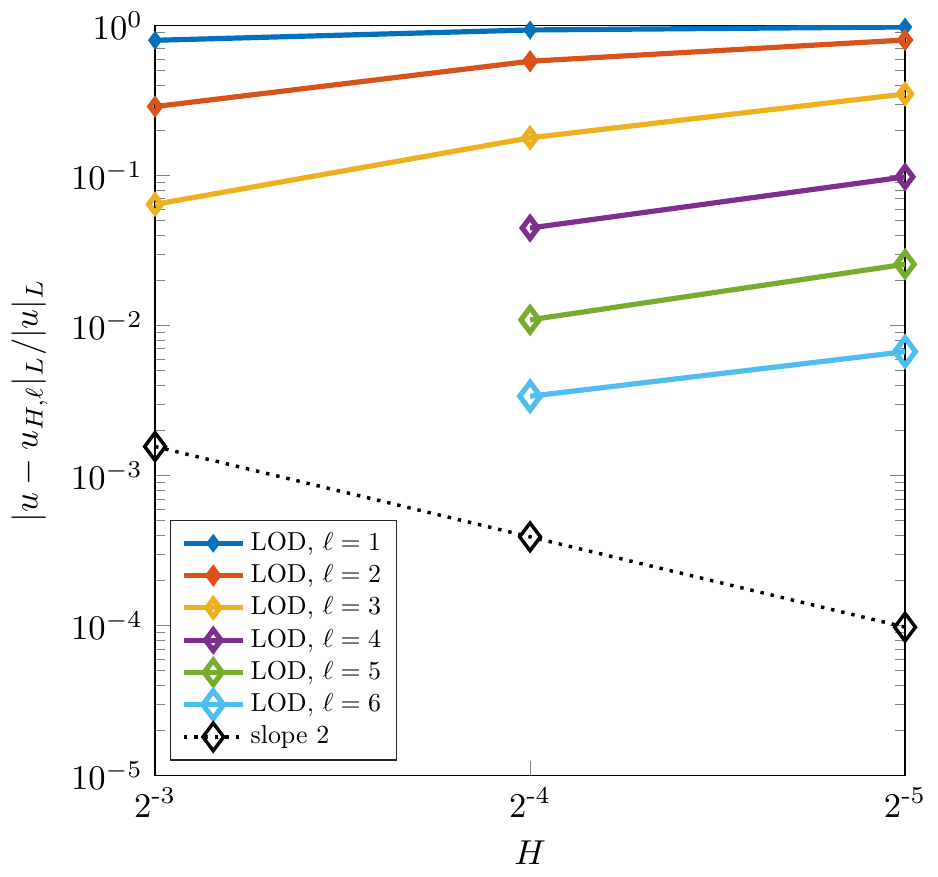}
	\caption{Plot of the errors $\lnormf{u - u_{H,\ell}}/\lnormf{u}$ in dependence of the mesh size $H$ of the SLOD (left) and the LOD (right).}
	\label{fig3}
\end{figure}

For demonstrating the convergence of the SLOD for spatial networks, we consider the edge weights from the previous numerical example and the right-hand side
\begin{equation*}
f(x_1,x_2) = \sin(x_1)\sin(x_2).
\end{equation*}
\Cref{fig3} shows the convergence plots for the SLOD (left) and the LOD (right) for different oversampling parameters as the coarse mesh $\TH$ is refined. Note that we only consider combinations of $H$, $\ell$ for which there is no patch $\mathsf N^\ell(T)$ that coincides with the whole domain~$\Omega$. As reference, lines of slope 2 are depicted.

\Cref{fig3} confirms the method's convergence properties as stated in \Cref{t:error}, i.e., provided that the oversampling parameter $\ell$ is chosen sufficiently large, optimal convergence of order 2 can be observed. Note that in \Cref{fig3} (left), the yellow line is partially below the purple line and thus is difficult to see. For the LOD, the optimal order convergence is difficult to see. The reason is that, for  all considered localization parameters, the localization error still dominates the optimal order error.

\subsection{High-contrast example}

\begin{figure}[h]
	\includegraphics[width=.42\linewidth]{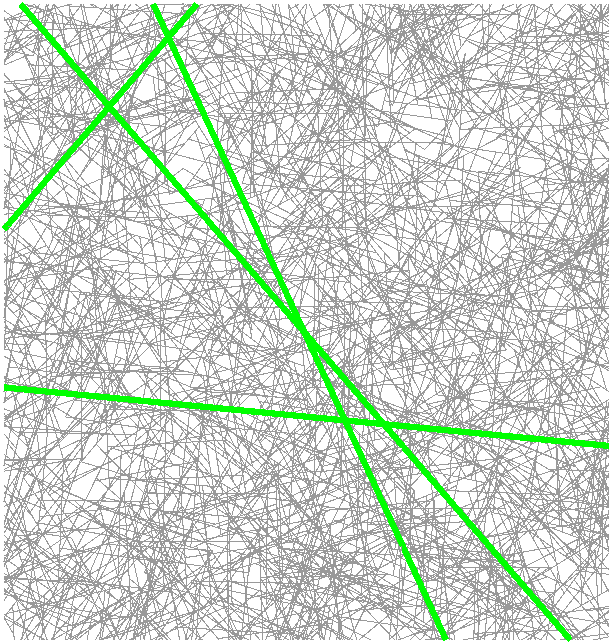}\hspace{1.5cm}
	\includegraphics[width=.42\linewidth]{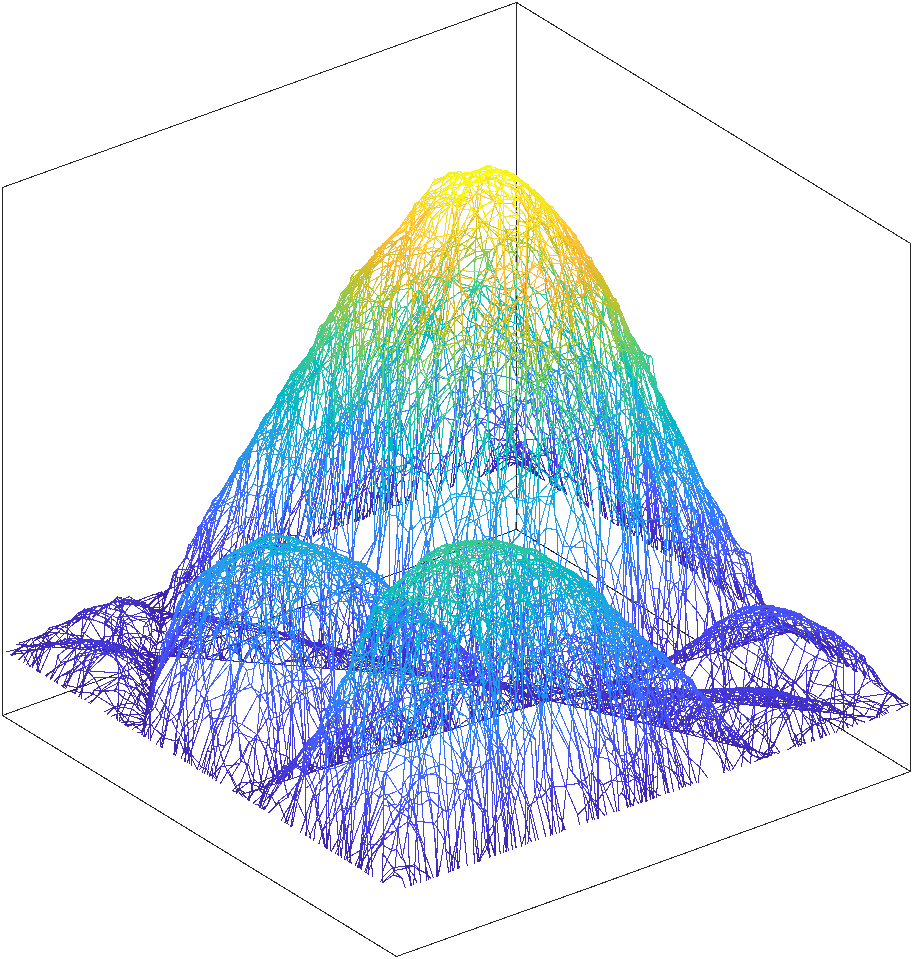}
	\caption{Spatial network with high conductivity channels indicated in green (left) and SLOD solution (right).}
	\label{fig:hc}
\end{figure}

In this numerical experiment, we use edge weights~$\gamma_{xy}$ that are independently and uniformly distributed in $[0.01,1]$ and add several  channels of high conductivity with edge weights of $10^4$. Hence, the contrast in this numerical example is of order $10^6$. For an illustration of the setup, we refer to \Cref{fig:hc} (left), where the high-conductivity channels are marked in green. 
The high conductivity effectively extends the homogeneous Dirichlet boundary conditions also to the channels. We consider this experiment as particularly challenging as the channels are not aligned with  the considered Cartesian meshes $\TH$ and thus, artificial boundary-like conditions are imposed on the SLOD basis functions. As right-hand side, we choose $f \equiv 1$, i.e., only the localization error is present, cf. \Cref{rem:exactness}.

For this experiment, the SLOD is again able to achieve very accurate approximations for much smaller oversampling parameters than the LOD. For example, for the mesh $\T_{2^{-4}}$ and $\ell = 3$, the SLOD achieves a relative  $L$-norm error of $2.75\times10^{-3}$, while the LOD, for the same discretization parameters, only achieves an error of $1.73 \times 10^{-1}$. For reaching a similar accuracy as the SLOD, for the LOD, we would need to choose  $\ell$ so large that many patch problems are already global problems.

We shall note that in \Cref{fig:hc}, for the ease of illustration, we only depicted a subnetwork and the correspondingly restricted solution. Plotting the full network is infeasible as the high density of the lines would make the lines nearly impossible to distinguish.

\section*{Acknowledgments}
We would like to thank Christoph Zimmer for inspiring discussion on the localization problem and especially its implementation.

\bibliographystyle{alpha}
\bibliography{bib}
\end{document}